\documentstyle[amsfonts,amssymb,graphicx]{article}  
 
\newcommand{\sgn}{\mathop{\mathrm{sgn}}}

\newtheorem{theorem}{Theorem}[section]
\newtheorem{lemma}[theorem]{Lemma}

\newtheorem{remark}[theorem]{Remark} 
\newtheorem{proposition}[theorem]{Proposition}

\title{\bf Global existence of the  self-interacting scalar field in  the de~Sitter universe}

\author{{\bf Karen Yagdjian}
 }

\begin{document}
\date{}
 
\maketitle
\thispagestyle{empty}
\vspace{-0.3cm}

\begin{center}
{\it School of Mathematical and Statistical Sciences,\\
University of Texas RGV,
1201 W.~University Drive, \\
Edinburg, TX 78539,
USA }
\end{center}
\medskip

\thispagestyle{empty} 
\vspace{-0.3cm}

\addtocounter{section}{-1}
\renewcommand{\theequation}{\thesection.\arabic{equation}}
\setcounter{equation}{0}
\pagenumbering{arabic}
\setcounter{page}{1}
\thispagestyle{empty}

\begin{abstract}
\begin{small}

We present some sufficient conditions for  the global in time existence of solutions of the semilinear Klein-Gordon equation  of the self-interacting  scalar field with complex mass. 
The coefficients of the equation depend on spatial variables as well, that makes results applicable, in particular, to the 
spacetime with the time slices being Riemannian manifolds.  The  least  lifespan  estimate   is given for the class of equations including the Higgs  boson equation,   which according to physics has a finite lifetime.  

\smallskip


\end{small}
\end{abstract}

\section{Introduction and Statement of Results}

\setcounter{equation}{0}
\renewcommand{\theequation}{\thesection.\arabic{equation}}

\smallskip

In this  paper we present some sufficient conditions for  the global in time existence of solutions of the semilinear Klein-Gordon equation  for the self-interacting  scalar field with complex mass. The estimate for the lifespan is given for the equation with the Higgs  potential. 
The coefficients of the equation depend on spatial variables as well, that makes results applicable, in particular, to the 
spacetime with the time slices being Riemannian manifolds. The case of the equation in the de~Sitter spacetime  (see, e.g., \cite[p.113]{Choquet-Bruhat_book}) is included.  
\smallskip

We consider the equation
\begin{eqnarray}
\label{0.1} 
&  &
  \psi  _{tt} +   n   \psi  _t - e^{-2 t} A(x,\partial_x) \psi   +  m^2  \psi  =   F(x,\psi  )\,,
 \end{eqnarray} 
where $ A(x,\partial_x) = \sum_{|\alpha |\leq 2} a_{\alpha }(x)\partial _x^\alpha $ is a  second order negative uniformly elliptic  operator with coefficients $a_{\alpha } \in {\mathcal B}^\infty $, where   ${\mathcal B}^\infty $  is the space of all $C^\infty ({\mathbb R}^n)$ functions with uniformly bounded derivatives of all orders.     We assume that the mass $m $ can be a complex number,   $m^2 \in{\mathbb C} $.  
\smallskip

In the quantum field theory    the description  of  matter fields is based on 
 the  semilinear Klein-Gordon equation generated by the mass $m$ and the metric $g$:
\begin{eqnarray*}
\square_g \psi = m^2 \psi   +   V'_\psi (x,\psi  ) \,.
\end{eqnarray*}
Here  $ \square_g $ is the Laplace-Beltrami operator. In physical
terms this equation describes a local self-interaction for a scalar particle.  
The special case of the equation (\ref{0.1}) is the   covariant Klein-Gordon equation   in the de~Sitter spacetime    
\[
\psi _{tt}
-  \frac{e^{-2t}}{\sqrt{|\det \sigma ( x)| }} \sum_{i,j=1}^n  \frac{\partial  }{\partial x^i}\left(  \sqrt{|\det \sigma  ( x)| }  \sigma  ^{ij} (x)\frac{\partial  }{\partial x^j}  \psi \right)
+ n   \psi_t +  m^2 \psi
 =
F(\psi ) \,.
\] 
The metric $\sigma    (x) $ belongs to the time slices. The metric $g$ in the de~Sitter   
spacetime   is  as follows, 
$g_{00}=  g^{00}= -  1  $, $g_{0j}= g^{0j}= 0$, $g_{ij}(x,t)=e^{2t}    \sigma  _{ij} (x) $, 
$i,j=1,2,\ldots,n$, where $\sum_{j=1}^n\sigma  ^{ij} (x) \sigma _{jk} (x)=\delta _{ik} $,   and $\delta _{ij} $ is  Kronecker's delta. 
\smallskip

 In \cite{JMAA_2012}-\cite{Helsinki_2013}  the global existence of  
 small data solutions of the Cauchy problem for the semilinear  Klein-Gordon equation and systems of equations  in the de~Sitter spacetime  with flat time slices,
was proved.
 The nonlinearity $ F$ was assumed Lipschitz continuous with the exponent $\alpha > 0 $ (see definition below) while  
   $m \in (0, \sqrt{n^2-1}/2]\cup [n/2,\infty) $.
The proof of the global existence in \cite{JMAA_2012}-\cite{Helsinki_2013} is based on the special integral representations (see Section~\ref{S2}) and $ L^p-L^q$ estimates. Later on,  in \cite{NARWA} this result for the same range of the  parameters $n,m$ and the same nonlinearity was extended on the equation (\ref{0.1}), that is,  from the spatially  flat de~Sitter spacetime to the    de~Sitter spacetime
with the time slices being, in particular, the Riemannian manifolds. The case of $m \in (\sqrt{n^2-1}/2,n/2) $  was left open in \cite{NARWA}.  The existence of solution in the energy spaces was not proved in \cite{NARWA}. Another interesting and important case, that is, the  case  of the  complex-valued mass $m$ also was not discussed in \cite{NARWA}. That case contains the  Klein-Gordon model of the Higgs boson equation.   
\smallskip

In the present paper we  generalize and complete  the  small data global existence result of \cite{NARWA}. In particular, we study also  class of equations containing the  Higgs boson equation with the Higgs potential, that is the equation
\begin{equation}
\label{Higgs}
\psi _{tt} - e^{-2t} A(x,D)\psi  +n\psi_t  =  \mu^2  \psi  - \lambda \psi  ^3 ,
\end{equation} 
 with $\lambda >0 $ and $\mu >0 $, while $n=3$. (For the Minkowski spacetime see, e.g., \cite[Ch.17]{Elbaz}.)
 
\smallskip

The  explicit form of  nonlinear term $F$ in this paper is    not used. What we use are simply 
the estimates of the form $\|F(\psi )\|_{X}   < C \| \psi  \|_{X'}^{\alpha }\| \psi  \|_{X'{}'}$, for some function
spaces $X$, $X'$ and $X'{}'$. 
Furthermore, since we prove the results for small data  in  the Sobolev space $H_{(s)}({\mathbb R}^n)$, we are only
concerned with the behavior of $F $ at the origin.  \\

\noindent {\bf Condition ($\mathcal L$).} {\it The smooth in $x$ function $F=F(x,\psi )$ is
said to be Lipschitz continuous with exponent $\alpha \geq 0 $ in the space $H_{(s)}({\mathbb R}^n)$   if there is a constant \,
$C \geq 0$ \, such that}
\begin{equation}
\label{calM}
 \|  F(x,\psi  _1 (x))- F(x,\psi  _2(x) ) \|_{H_{(s)} }  \leq C
\| \psi  _1 -  \psi _2   \|_{H_{(s)} }
\Big( \|  \psi _1  \|^{\alpha} _{H_{(s)} }
+ \|  \psi _2   \|^{\alpha} _{H_{(s)} } \Big) \,\, \mbox{\it for all} \,\,  \psi  _1, \psi  _2 \in  H_{(s)}\,.
\end{equation} 
\smallskip

\noindent
The  polynomials  in $ \psi$ are Lipschitz continuous with some exponent $\alpha $ in the space $H_{(s)}({\mathbb R}^n)$ when $s>n/2 $. Moreover, the exponent $\alpha $ is independent of $s$.  
Interesting  functions  are   $F(x,\psi )=\pm|\psi |^{\alpha +1}$,  $F(\psi )=\pm|\psi  |^{\alpha } \psi $
as important examples of  Lipschitz continuous functions in the Sobolev space $H_{(s)}({\mathbb R}^n) $ for $\alpha > 0 $, $s>n/2$,  provided that $\alpha $  agrees with $s $ and $n$. More detailed interplay between $\alpha $, $s $, and $n$  of the Condition   ($\mathcal L$) is an issue interesting in its own right   but it is out of the scope of this paper. 

Define also the metric space
\[
X({R,H_{(s)},\gamma})  := \left\{ \psi  \in C([0,\infty) ; H_{(s)} ) \; \Big|  \;
 \parallel  \psi   \parallel _X := \sup_{t \in [0,\infty) } e^{\gamma t}  \parallel
\psi  (x ,t) \parallel _{H_{(s)}}
\le R \right\}\,,
\]
where $\gamma \in {\mathbb R}$, with the metric
\[
d(\psi _1,\psi _2) := \sup_{t \in [0,\infty) }  e^{\gamma t}  \parallel  \psi _1 (x , t) - \psi _2 (x ,t) \parallel _{H_{(s)}}\,.
\]

We study the Cauchy problem   (\ref{NWE}), (\ref{ICPHI})
    through the integral equation.
To define that integral equation we  appeal to the operator 
\[
G:={\mathcal K}\circ {\mathcal EE}  
\]
(${\mathcal EE}$ stands for the evolution equation) as follows. For the function $f(x,t) $ we define
\[
v(x,t;b):= {\mathcal EE} [f](x,t;b)\,,
\]
where the function 
$v(x,t;b)$   
is a solution to the Cauchy problem 
\begin{eqnarray}
\label{1.6} 
&   &
\partial_t^2 v - A(x,D)v =0, \quad x \in {\mathbb R}^n, \quad t \geq 0, \\
\label{2.2a}
&  &
v(x,0;b)=f(x,b)\,, \quad v_t(x,0;b)= 0\,, \quad x \in {\mathbb R}^n\,,
\end{eqnarray} 
while ${\mathcal K}$ is introduced  by 
\begin{eqnarray}
\label{0.5} 
{\mathcal K}[v]  (x,t) 
&  :=  &
2   e^{-\frac{n}{2}t}\int_{ 0}^{t} db
  \int_{ 0}^{ e^{-b}- e^{-t}} dr  \,  e^{\frac{n}{2}b} v(x,r ;b) E(r,t; 0,b;M)  \,.
\end{eqnarray}
The kernel $ E(r,t; 0,b;M) $ was introduced in \cite{JMAA_2012} and \cite{Yag_Galst_CMP} (see also (\ref{E})). Hence, 
\[
G[f]  (x,t)  
  =  
2   e^{-\frac{n}{2}t}\int_{ 0}^{t} db
  \int_{ 0}^{ e^{-b}- e^{-t}} dr  \,  e^{\frac{n}{2}b}\,{\mathcal EE} [f](x,r ;b) E(r,t; 0,b;M)  \,.
\] 
Thus, the Cauchy problem (\ref{NWE}), (\ref{ICPHI}) leads to the following integral equation
\begin{eqnarray} 
\label{5.1}
\psi  (x,t)
 = 
\psi  _0(x,t) + 
G[ F(\cdot ,\psi ) ] (x,t)    \,. 
\end{eqnarray}
Every solution to  the  Cauchy problem  (\ref{NWE})-(\ref{ICPHI}) solves also the last integral equation with some function $\psi  _0 (x,t)$, which is a solution for the problem for the linear equation without source term. We define a solution of the Cauchy problem (\ref{NWE})-(\ref{ICPHI})  via integral equation (\ref{5.1}). 
Since only for $m \in(0,\sqrt{n^2-1}/2]\cup [n/2,\infty)  $ the existence of global in time solution has been proved in \cite{NARWA}, in the present paper we  consider the more general case of  the complex mass $m \in {\mathbb C}$ that includes, in particular, the Higgs boson equation.
The principal square root $M:=(n^2/4-m^2)^\frac{1}{2} $ is the  parameter that controls estimates and solvability. In fact, 
${\mathcal M} :=iM $ is the so-called {\it effective mass} or {\it curved mass} of the field. The main result of this paper is the next theorem. 
\begin{theorem}
\label{T0.1}
Assume that the nonlinear term $F(x,\psi  )$  is a  Lipschitz continuous in the  space $H_{(s)} ({\mathbb R}^n)$, $ s > n/2\geq 1$, $F(x,0)\equiv 0$, and  $\alpha >0 $.

\noindent
$\mbox{\rm (i)}$ Assume also that   $0< \Re M <1/2$.   
Then, there exists $\varepsilon _0>0 $ such that, for every given functions $\psi _0 ,\psi _1 \in H_{(s)} ({\mathbb R}^n) $,  such that 
\begin{equation}
\label{0.7}
 \| \psi _0   \|_{ { H}_{(s)} ({\mathbb R}^n)}
+  \|\psi _1  \|_{ { H}_{(s)} ({\mathbb R}^n)} \leq \varepsilon, \qquad \varepsilon  < \varepsilon_0\,, 
\end{equation} 
there exists a solution $\psi  \in C ([0,\infty);H_{(s)} ({\mathbb R}^n))$ of the Cauchy problem  
\begin{eqnarray}
\label{NWE}  
&  &
  \psi  _{tt} +   n   \psi  _t - e^{-2 t} A(x,\partial_x)  \psi   +  m^2  \psi  =   F(x, \psi  )\,,\\
\label{ICPHI} 
&  &
 \psi  (x,0)= \psi _0 (x)\, , \quad \psi   _t(x,0)=\psi _1 (x) \,.
 \end{eqnarray} 
The  solution \, $ \psi   (x ,t) $ \, belongs to the space \, $  X({2\varepsilon,s, \frac{n-1}{2}  })  $, that is,  
\begin{eqnarray*}  
\sup_{t \in [0,\infty)}  e^{ \frac{n-1}{2} t}  \|\psi  (\cdot ,t) \|_{H_{(s)} ({\mathbb R}^n)}  \leq  2\varepsilon\, .
\end{eqnarray*}

\noindent
$\mbox{\rm (ii)}$ Assume  that  $M=1/2$ or  $1/2 <  \Re M <n/2$ and $  \gamma \in (0,\frac{1}{\alpha +1}(\frac{n}{2}- \Re M ))$. Then there exists $\varepsilon _0>0 $ such that  for every given functions $\psi _0 ,\psi _1 \in H_{(s)} ({\mathbb R}^n) $,  satisfying (\ref{0.7}), 
there exists a solution $\psi  \in X({2\varepsilon,s, \gamma })$ of the Cauchy problem (\ref{NWE})-(\ref{ICPHI}).\\ 

\noindent
$\mbox{\rm (iii)}$ If $\Re M > n/2$, then  
the lifespan $T_{ls}$ of the solution can be estimated from below as follows
\begin{eqnarray*} 
   T_{ls}    
&  \geq   &
-\frac{1}{ \Re M  - \frac{n}{2} }\ln  \left( \| \psi _0   \|_{ { H}_{(s)} ({\mathbb R}^n)}
+  \|\psi _1  \|_{ { H}_{(s)} ({\mathbb R}^n)} \right)     
-C(m,n,\alpha ) 
\end{eqnarray*} 
with some constant $C(m,n,\alpha ) $.
 \end{theorem}
In particular, the theorem covers the case of    
$m \in (\sqrt{n^2-1}/{2}, n /{2})$, which was left open in \cite{NARWA}. 
If 
\[
F(\psi )=\lambda \psi  ^3
\quad \mbox{\rm or}  \quad F(\psi )= \pm  |\psi  |^\alpha \psi  \quad \mbox{\rm or}  \quad F(\psi )= \pm  |\psi |^{\alpha+1},
\]
then the small data Cauchy problem is  globally  solvable for every $\alpha  $, $s $, and $n$ satisfying ($\mathcal L$).
 \medskip

We note that there is some discontinuity  at $M=1/2 $  in the decay rate of the solution in the transition from the  cases (i) to (ii). It is a result  of the nonlinearity since if  $\alpha  =0$, then the discontinuity  disappears. On the other hand    
the invariance under some  gauge  transformation  already suggests that there exists some discontinuity in the theory at $m =\sqrt{2}$ if $n=3$, that is, if $M=1/2 $ (see \cite{Deser-Waldron}).  
\medskip
 
The   finite time blowup of the solutions is proved in \cite{yagdjian_DCDS} for the Cauchy problem (\ref{NWE})-(\ref{ICPHI}) with the  wide class of semilinear term $F$. That class is contained in one assumed in (iii) of Theorem~\ref{T0.1}. In particular, the blow up occurs  if $M>n/2$ and $F(x,\psi )=|\psi |^\alpha $ with $\alpha >0$.   The last class of $F$ includes  sign preserving solutions of the 
Higgs boson equation (\ref{Higgs}). 
In fact, the application of the  assertion (iii) of Theorem~\ref{T0.1} to the Higgs boson equation 
 gives an estimate from below of the lifespan of the Higgs boson. 
The mathematical proof of the finiteness of the lifespan of solution, which is not necessarily sign preserving, to the equation (\ref{Higgs}) is an interesting and difficult problem that  requires special technique
(see, e.g., \cite{yagdjian_DCDS})  and   will be  published in a forthcoming  paper. In the existing extensive literature on the physics of the Higgs boson one can find that   the lifetime of the Higgs boson is approximately $10^{-13}sec.$ (see, e.g., \cite{CMS_Collaboration,CMS_Collaboration2}) that points at the boundedness of the  lifespan of some solutions to the equation (\ref{Higgs}).  
\medskip

Although, there is no conservation of energy due to the dependence on time of the coefficient, 
 the energy estimate provides with the useful tool to prove global existence in the energy space if we impose some restriction on the nonlinearity. The last theorem as well as the results of articles \cite{Galstian-Yagdjian-NA,NARWA} imply global solvability of the problem  in the energy space under some conditions on the nonlinear term $F$ and mass $m$.  
\begin{theorem}
\label{T0.2}
Assume that the nonlinear term $F $  is   Lipschitz continuous in the  space $H_{(s)} ({\mathbb R}^n)$, $ s > n/2\geq 1$, $F(0)=0$, and  $\alpha >\frac{2}{n-1} $. 
Assume also that either $M^2 \in  {\mathbb R}$ and $ \Re M \in (0,1/2) $ or $M=1/2$.  
Then, there exists $\varepsilon _0>0 $ such that, for every given functions $\psi _0  \in H_{(s+1)} ({\mathbb R}^n) $, $\psi _1 \in H_{(s)} ({\mathbb R}^n) $, such that 
\[
 \| \psi _0   \|_{ { H}_{(s+1)} ({\mathbb R}^n)}
+  \|\psi _1  \|_{ { H}_{(s)} ({\mathbb R}^n)} \leq \varepsilon, \qquad \varepsilon  < \varepsilon_0\,, 
\] 
there exists a global solution $\psi  \in C^1([0,\infty);H_{(s)} ({\mathbb R}^n))$ of the Cauchy problem  
(\ref{NWE})-(\ref{ICPHI}). 
The  solution \, $ \psi   (x ,t) $ and its time derivative $ \partial_t\psi   (x ,t)$  belong to the space \, $  X({2\varepsilon,s, \frac{n-1}{2 }  })  $, that is,  
\begin{eqnarray*}  
\sup_{t \in [0,\infty)}  e^{ \frac{n-1}{2 }   t}  \left( \|\psi   (\cdot ,t) \|_{H_{(s)} ({\mathbb R}^n)} 
+\| \partial_t \psi  (\cdot ,t) \|_{H_{(s)} ({\mathbb R}^n)} 
\right) < 2\varepsilon \,.
\end{eqnarray*} 
Assume that  $M^2 \in  {\mathbb R}$ and either $ \Re M \in (3/2,n/2) $ or $M=3/2$, then there exists $\varepsilon _0>0 $ such that, for every given functions $\psi _0  \in H_{(s+1)} ({\mathbb R}^n) $, $ \psi _1 \in H_{(s+1)} ({\mathbb R}^n) $, such that 
\[
 \|\psi _0   \|_{ { H}_{(s+1)} ({\mathbb R}^n)}
+  \|\psi _1  \|_{ { H}_{(s+1)} ({\mathbb R}^n)} \leq \varepsilon, \qquad \varepsilon  < \varepsilon_0\,, 
\] 
there exists a global solution $\psi \in C^1([0,\infty);H_{(s)} ({\mathbb R}^n))$ of the Cauchy problem  
(\ref{NWE})-(\ref{ICPHI}) such that  \, $ \psi   (x ,t)\in X({2\varepsilon,s, \gamma   })  $ and its time derivative $ \partial_t\psi   (x ,t)$  belong to the space \, $  X({2\varepsilon,s, \gamma    -1})  $ with $\gamma \in (0, \frac{1}{\alpha +1}(\frac{n}{2}-\Re M )) $, that is,  
\begin{eqnarray*}  
\sup_{t \in [0,\infty)}  e^{  \gamma   t}     \|\psi   (\cdot ,t) \|_{H_{(s)} ({\mathbb R}^n)} 
+\sup_{t \in [0,\infty)} e^{  (\gamma-1) t} \| \partial_t \psi   (\cdot ,t) \|_{H_{(s)} ({\mathbb R}^n)} 
  < 2\varepsilon \,.
\end{eqnarray*}

 \end{theorem}
\smallskip

The main tools to study the problem (\ref{NWE})-(\ref{ICPHI}) are the integral transform from \cite{MN} and the standard energy estimate for the finite time interval for the strictly hyperbolic equation. On the other hand, by using   the integral transforms given in \cite{MN}, it is possible to reduce    the problem with the infinite time interval to the problem for the hyperbolic equation with time independent coefficients and with the finite time  interval due to the fact that    the de~Sitter spacetime has permanently bounded domain of influence. In this approach  the integral transform allows us to push forward the  estimates provided that some integrals of the kernel functions lead to the proper estimates. The proof of the estimates for the integrals of kernel the  functions consists of a long sequence of estimates of integrals involving hypergeometric functions. The proof of Theorem~\ref{T0.1} is concluded by   
fixed point arguments.

\smallskip

The  estimates derived for the linear equation in Sections~\ref{S2}, \ref{S3} include equations of scalar fields considered in \cite{Bros-2,Epstein-Moschella} with $m^2<0 $  living on the
de~Sitter universe. 
The Klein-Gordon scalar quantum fields on the de~Sitter manifold  with 
imaginary mass   $m^2=-k(k+n)$, $\,k=0,1,2,\ldots$, 
present a family of tachyonic  quantum fields.    Epstein and Moschella \cite{Epstein-Moschella} give a complete study of such family of linear scalar tachyonic quantum
fields.   The corresponding linear equation is
\begin{eqnarray}
\label{0.10new}
&  &
\psi _{tt} +   n  \psi  _t - e^{-2 t} \Delta  \psi   +  m^2  \psi  =  0\,,  
 \end{eqnarray} 
for which the kernel of the integral transform ${\mathcal K}$ (\ref{0.5}) is  $E(x,t;x_0,t_0;M) $, where  $M  = k+\frac{n }{2}$, $k=0,1,2,\ldots\,$. For  an odd number $n$,  the mass 
$m$ takes values on the  set of  knot points in the sense of \cite{JMP2013}.    Theorem~\ref{T0.1} contains an estimate for  the lifespan of a self-interacting tachyonic field. 
\medskip

The approach  of this paper can be easily modified to obtain estimates for the linear equation and the  global (in time) solvability for the equation 
\begin{eqnarray} 
\label{0.1b}
&  &
 \psi _{tt} +   \nu   \psi _t - e^{-2 t} A(x,\partial_x) \psi  +  m^2 \psi  = f(x,t)+  F(x,\psi )  \,,
 \end{eqnarray} 
 where $\nu \in {\mathbb R} $. It is also possible include  the derivatives of the field function in the nonlinear term,  
as well as to apply the approach  of the present paper to system of equations similar to \cite{Helsinki_2013}.

\smallskip

The equation (\ref{0.1}) in the case of  the  $x$-independent operator $A(x,\partial_x)=\Delta  $ is amenable to the analysis via the Fourier transform and  the Bessel functions (see, e.g,  \cite{Galstian}). For the equation (\ref{0.1b}) with $A(x,\partial_x)=\Delta  $, $f(x,t)=0$, and $F(x,\psi )= |\psi  |^p$ the Fourier transform was used  in \cite{E-R}, and global existence in the energy classes and in the  Lebesgue spaces was proved under several restrictions on the nonlinear term. The $x$-independence of the coefficients
allows authors to apply the Fourier transform and to write an explicit form of the solution of the corresponding ordinary differential equation.       
\smallskip

Unlike to  the  case of the  operator $A(x,\partial_x)=\Delta  $, the linear part of the equation (\ref{0.1})  is not invariant with respect to de~Sitter group (see, e.g., \cite{Garidi,GOTO}). Nevertheless, the mass intervals $(0, \sqrt{n^2-1}/2)$, $[\sqrt{n^2-1}/2,n/2)$, $[n/2,\infty) $ appear and play important role also in this case. The first interval $(0, \sqrt{n^2-1}/2)$ with $n=3$  in quantum field theory is known as  the  Higuchi bound  ({\it forbidden mass range}). 
The masses in this range  lead to negative norm states, i.e., non-unitarity. In \cite{Higuchi} it is shown that for spin-2 fields the forbidden mass range is $0<m^2<2$. The mass $m= \sqrt{n^2-1}/2 $ is remarkable especially because that is the only mass that makes equation (\ref{0.10new}) Huygensian \cite{JMP2013} and makes the linear part of the equation conformally  invariant \cite{Birrell}. The values   $ 0$ and $\sqrt{n^2-1}/2  $ are the only values of mass such that the equation obeys {\it incomplete Huygens' principle} \cite{JMP2013}. 
In the de~Sitter spacetime the existence of two different scalar fields (in fact, with $m=0$ and $m^2=(n^2-1)/4 $), 
which obey  incomplete Huygens' principle, 
is equivalent  to the condition $n=3$ (Corollary~4~\cite{JMP2013}), which is the spatial dimension of the physical world. 
In fact, Paul Ehrenfest
in \cite{Ehrenfest} addressed the question: ``Why has our space just three dimensions?''. 
\smallskip
 
Thus, the point $m= \sqrt{2} $ ($n=3$)  is exceptional for the  quantum fields theory in the de~Sitter spacetime.  In particular, for massive spin-2 fields, it is known \cite{Deser-Waldron,Higuchi} that the norm of the helicity zero mode changes sign across the line $m^2= 2 $. The region  $m^2<2$ is therefore unitarily forbidden. It is noted in \cite{Lasma Alberte} that all canonically normalized helicity $-0,\pm 1,\pm 2\,$ modes of massive graviton on the de Sitter universe satisfy Klein-Gordon equation for a {\it massive scalar field with the same effective mass}.   For the case of large mass, that is $m^2\geq n^2/4 $, and for the brief review of the bibliography related to that case, one can consult \cite{NARWA,Nakamura,Spr2016} and for the results on the equation in the asymptotically de~Sitter spaces see \cite{Baskin,BaskinSE,Hintz-Vasy,Hintz,Vasy_2010}. The waves in spacetimes with a nonvanishing cosmological constant are studied in \cite{Bachelot_2016,Costa,Jetzer}.
\smallskip

If $n=3$, then  another important value  is $m= 3/2$.  
The equation for the scalar field with mass  $m$   
in de~Sitter universe in the physical variables is:    
\[
\frac{1}{c^2}\psi _{tt} +   \frac{1}{c^2} 3 H   \psi _t - e^{-2tH} \bigtriangleup \psi  + \left( \frac{c m}{h}\right)^2\psi = 0\,.  
\]  
Here $h = 1.054\cdot 10^{-27} erg \cdot sec$, $c\approx 3\cdot 10^{10} \frac{cm}{sec}$, $H \approx 10^{-18} \frac{1}{sec}$. 
The following question seems to be natural:  For what particle (mass) the equation has the most simple form? 
In fact, tor the scalar field with the mass  $m= \frac{3hH}{2c^2}$ the function $u = e^{-\frac{3}{2}Ht}\psi $  solves the equation  
\[
\frac{1}{c^2}u_{tt}     - e^{-2tH} \bigtriangleup u =  0\,. 
\] 
In the physical units this particle has a mass
$
m= \frac{3hH}{2c^2}\approx 1.756 \cdot 10^{-66} 
g$.
The natural   question  arises: What particle has this mass?  
In fact, there exists an extensive literature on this topic.  
The comparison  with the estimate   $m_g < 1.8\cdot 10^{-66}g$  from \cite{Genk-Tron} supports the following conjecture (see, e.g., \cite{yagdjian_DCDS}): the mass  
$\displaystyle m= 3hH/(2c^2)$  is a mass of graviton.  
\smallskip

The present paper is organized as follows. In Section~\ref{S2} we describe 
 the  integral transform and the   generated by that transform  representations (from \cite{MN})
for the solutions of the Cauchy problem for the linear equation. Then, 
we  show that the energy estimates for the second order hyperbolic operator with time independent coefficients
can be pushed forward via integral transform 
to the source free equation with time dependent coefficient.  
In  the present paper we prove estimates for the Sobolev spaces only. In fact, the proofs for the Lebesgue, Sobolev and Besov spaces are identical.   In Section~\ref{S2b} we  obtain similar estimates for the equation with  source term. 
  The last section, Section~\ref{S3}, is devoted to the solvability of the associated integral equation and to the proof of Theorem~\ref{T0.1} and 
 Theorem~\ref{T0.2}.
In the Appendix one can find several useful lemmas concerning hypergeometric functions which have been used in the previous sections.

\section{$H_{(s)} ({\mathbb R}^n)$ Estimates}
\label{S2}

\setcounter{equation}{0}

We introduce  the kernel functions $E(x,t;x_0,t_0;M) $, $K_0(z,t;M)   $,    and $K_1(z,t;M) $    (see also \cite{Yag_Galst_CMP} and \cite{JMAA_2012}). 
 First, for $M \in {\mathbb C} $  we define the function 
\begin{eqnarray}
\label{E}  
E(x,t;x_0,t_0;M)
 &  = & 
 4 ^{-M}  e^{ M(t_0+t) } \Big((e^{-t }+e^{-t_0})^2 - (x - x_0)^2\Big)^{-\frac{1}{2}+M    } \\
 &  &
 \times 
F\Big(\frac{1}{2}-M   ,\frac{1}{2}-M  ;1; 
\frac{ ( e^{-t_0}-e^{-t })^2 -(x- x_0 )^2 }{( e^{-t_0}+e^{-t })^2 -(x- x_0 )^2 } \Big) . \nonumber 
\end{eqnarray} 
Next 
we define also the kernels  $K_0(z,t;M)   $    and $K_1(z,t;M) $ by
\begin{eqnarray*}
K_0(z,t;M)
& := & 
- \left[  \frac{\partial }{\partial b}   E(z,t;0,b;M) \right]_{b=0} \\
&  = &   
4 ^ {-M}  e^{ t M}\big((1+e^{-t })^2 - z^2\big)^{  M    } \frac{1}{ [(1-e^{ -t} )^2 -  z^2]\sqrt{(1+e^{-t } )^2 - z^2} }\\
&   &
\times  \Bigg[  \big(  e^{-t} -1 +M(e^{ -2t} -      1 -  z^2) \big) 
F \Big(\frac{1}{2}-M   ,\frac{1}{2}-M  ;1; \frac{ ( 1-e^{-t })^2 -z^2 }{( 1+e^{-t })^2 -z^2 }\Big) \\
&  &
\hspace{1cm}  +   \big( 1-e^{-2 t}+  z^2 \big)\Big( \frac{1}{2}+M\Big)
F \Big(-\frac{1}{2}-M   ,\frac{1}{2}-M  ;1; \frac{ ( 1-e^{-t })^2 -z^2 }{( 1+e^{-t })^2 -z^2 }\Big) \Bigg]
\end{eqnarray*} 
and $K_1(z,t;M)   :=  
  E(z ,t;0,0;M) $, that is, 
\begin{eqnarray*} 
K_1(z,t;M)  
& = &
  4 ^{-M} e^{ Mt }  \big((1+e^{-t })^2 -   z  ^2\big)^{-\frac{1}{2}+M    } \\
  &  &
  \times 
F\left(\frac{1}{2}-M   ,\frac{1}{2}-M  ;1; 
\frac{ ( 1-e^{-t })^2 -z^2 }{( 1+e^{-t })^2 -z^2 } \right), \, 0\leq z\leq  1-e^{-t}, 
 \end{eqnarray*}  
respectively.  
The solution $\psi  $ to the Cauchy problem  
\begin{equation}
\label{1.3}
  \psi  _{tt} +   n   \psi  _t - e^{-2 t} A(x,\partial_x) \psi   + m^2\psi  =  f ,\quad \psi  (x,0)= \psi _0(x)  , \quad \psi _t(x,0)=\psi _1(x), 
\end{equation}
with \, $ f \in C^\infty ({\mathbb R}^{n+1})$\, and with  \, $ \psi _0 $,  $ \psi _1 \in C_0^\infty ({\mathbb R}^n) $, $n\geq 2$,  is given in \cite{MN} by the next expression 
\begin{eqnarray}
\label{Phy}
\psi  (x,t) 
 & =  &
2   e^{-\frac{n}{2}t}\int_{ 0}^{t} db
  \int_{ 0}^{ e^{-b}- e^{-t}} dr  \,  e^{\frac{n}{2}b} v(x,r ;b) E(r,t; 0,b;M)  \\
&  &
+ e^{-\frac{n-1}{2}t} v_{\psi _0}  (x, \phi (t))
+ \,  e^{-\frac{n}{2}t}\int_{ 0}^{1} v_{\psi _0}  (x, \phi (t)s)\big(2  K_0(\phi (t)s,t;M)+ nK_1(\phi (t)s,t;M)\big)\phi (t)\,  ds  \nonumber \\
& &
+\, 2e^{-\frac{n}{2}t}\int_{0}^1   v_{\psi _1 } (x, \phi (t) s) 
  K_1(\phi (t)s,t;M) \phi (t)\, ds
, \quad x \in {\mathbb R}^n, \,\, t>0\, , \nonumber 
\end{eqnarray}
where the function 
$v(x,t;b)$   
is a solution to the Cauchy problem (\ref{1.6})-(\ref{2.2a}),   
while $\phi (t):=  1-e^{-t} $.  Here, for $\varphi \in C_0^\infty ({\mathbb R}^n)$ and for $x \in {\mathbb R}^n$, the function $v_\varphi  (x, \phi (t) s)$  coincides with the value $v(x, \phi (t) s) $ 
of the solution $v(x,t)$ of the Cauchy problem for the equation (\ref{1.6})  with the initial datum $\varphi  (x) $ while the second datum is zero.
\medskip

The mass $m^2= (n^2-1)/4 $, that is,  $M=1/2$,      
    simplifies the hypergeometric functions, as well as, the kernels  $E\left(x,t;x_0,t_0;\frac{1}{2}\right)$, $K_0 (z,t;M)$ and $K_1 (z,t;M)$ 
(see \cite{JMP2013}). In that case 
\begin{eqnarray*}
E\left(x,t;x_0,t_0;\frac{1}{2}\right) = 
 \frac{1}{2} e^{ \frac{1}{2}(t_0+t) } ,\quad K_0\left(z,t;\frac{1}{2} \right)
  =
- \frac{1}{4}  e^{ \frac{1}{2}t }  ,\qquad
K_1\left(z,t;\frac{1}{2} \right)
   =   \frac{1}{2}  e^{ \frac{1}{2}t }  \,.
\end{eqnarray*}
For the solution   of the Cauchy problem (\ref{1.3})   it follows
\begin{eqnarray}
\label{1.5}
\psi  (x,t)
&  =  & 
    e^{-\frac{n-1}{2}t}\int_{ 0}^{t}  e^{\frac{n+1}{2}b} db
  \int_{ 0}^{ e^{-b}- e^{-t}}v(x,r ;b)  \,   dr +e^{-\frac{n-1}{2}t} v_{\psi _0}  (x, 1-e^{-t})\\
  &  &
+ \,   \frac{n-1}{2}e^{-\frac{n-1}{2}t}\int_{ 0}^{1-e^{-t}} v_{\psi _0}  (x, s )        \,  ds  
+\,  e^{-\frac{n-1}{2}t}\int_{0}^{1-e^{-t}}   v_{\psi _1 } (x, s )
     \, ds
, \quad x \in {\mathbb R}^n, \,\, t>0  \nonumber \,,
\end{eqnarray}
where the functions $v(x,r ;b) $,  $v_{\varphi_0 }(x, s )     $, and $v_{\varphi _1 } (x, s )  $ are defined above. 
\smallskip

\subsection{$\bf H_{(s)} ({\mathbb R}^n)$ Estimates for  Equations without Source}

Let $ A(x,\partial_x) = \sum_{|\alpha |\leq 2} a_{\alpha }(x)\partial _x^\alpha $ be a  second order negative uniformly elliptic  operator with coefficients $a_{\alpha } \in {\mathcal B}^\infty $, where   ${\mathcal B}^\infty $  is the space of all $C^\infty ({\mathbb R}^n)$ functions with uniformly bounded derivatives of all orders.  
 Let $u=u(x,t) $ be the solution of 
\begin{eqnarray}
\label{2.1}
&   &
\partial_t^2 v - A(x,D)v =0, \quad x \in {\mathbb R}^n, \quad t \geq 0, \\
\label{2.2}
&  &
v(x,0) = v_0 (x),\quad   v_t (x,0)= v_1(x),\quad x \in {\mathbb R}^n\,.
\end{eqnarray}
 The following energy estimate is well known. (See, e.g., \cite{Taylor}.) For every $s \in {\mathbb R}$ there is  $ C_s$ such that  
\begin{eqnarray}
\label{Blplq}
&   &
 \|v_{t}(t)\|_{H_{(s)}}   +  \|v (t)\|_{H_{(s+1)}} \leq  C_s ( \|v_1\|_{H_{(s)}} +  \|v_0 \|_{H_{(s+1)}} )  , \quad 0 \leq  t \leq 1  \,.
\end{eqnarray}

We note that although in this estimate the time interval is bounded, meanwhile, due to the integral transforms given in \cite{MN}, it is possible to reduce    the problem with infinite time to the problem with the finite time, and to apply (\ref{Blplq}). 
We must  to emphasize  that this is possible since the de~Sitter spacetime has permanently bounded domain of influence. 
    
\begin{theorem} 
\label{T13.2}
For every given $s \in {\mathbb R}$, 
the solution  $ \psi = \psi (x,t)$ of the Cauchy problem 
\begin{equation}
\label{1.7}
  \psi _{tt} +   n    \psi _t - e^{-2 t} A(x,D)  \psi   + m^2 \psi  =  0\,, \quad  \psi  (x,0)= \psi _0 (x)\, , \quad  \psi _t(x,0)=\psi _1 (x)\,,
\end{equation} 
with  $\Re M= \Re (\frac{n^2}{4}-m^2)^{1/2}\in (0,1/2)$ 
satisfies the following  estimate
\[ 
\|  \psi (x,t) \| _{H_{(s)}}  
  \leq   
C_{m,n,s}  e^{   -\frac{n-1}{2} t}\Big\{   \| \psi _0   
\|_{H_{(s)}}
+ (1- e^{-t}) \|\psi _1  
\|_{H_{(s)}} 
 \Big\} \quad \mbox{for all} \quad t \in (0,\infty)\,.
\]

If   
$\Re M= \Re (\frac{n^2}{4}-m^2)^{1/2}>  1/2$ or $M=1/2$,  
then    the solution  $ \psi  = \psi (x,t)$ of the Cauchy problem (\ref{1.7}) 
satisfies the following  estimate
\begin{eqnarray*}   
\| \psi  (x,t) \| _{H_{(s)}} 
   \leq 
C  e^{(\Re M-\frac{n}{2})t}     \left\{ \|\psi _0 \|_{H_{(s)}} + (1- e^{-t}) \|\psi _1 \|_{H_{(s)}}\right\} \quad \mbox{for all} \quad t \in (0,\infty)\,.
\end{eqnarray*}  
\end{theorem}

\noindent
{\bf Proof.} The case of $M=1/2 $ is an evident consequence of (\ref{Blplq}) and the representation (\ref{1.5}) and in the remaining part of the proof it is not discussed. 

First 
we consider the case of $\psi _1=0 $. Then
\[
 \psi  (x,t) 
  =  
e^{-\frac{n-1}{2}t} v_{\psi _0}  (x, \phi (t))  
+ \,  e^{-\frac{n}{2}t}\int_{ 0}^{1} v_{\psi _0}  (x, \phi (t)s)\big(2  K_0(\phi (t)s,t;M)+ nK_1(\phi (t)s,t;M)\big)\phi (t)\,  ds
\]
and, consequently, 
\begin{eqnarray} 
\label{2.9}
\|  \psi  (x,t) \| _{H_{(s)}}  
&   \leq &  
e^{-\frac{n-1}{2}t} \|  v_{\psi _0}  (x, \phi (t)) \| _{H_{(s)}}   \\
&  &
+ \,  e^{-\frac{n}{2}t}\int_{ 0}^{1} \| v_{\psi _0}  (x, \phi (t)s)\| _{H_{(s)}} \big|2  K_0(\phi (t)s,t;M)+ nK_1(\phi (t)s,t;M)\big|\phi (t)\,  ds\,. \nonumber
\end{eqnarray}
Then  for the solution $v = v (x,t)$ of the Cauchy problem (\ref{2.1})-(\ref{2.2}) 
with $\varphi (x) \in C_0^\infty({\mathbb R}^n)$ one has  the estimate (\ref{Blplq}).
Hence, 
\begin{eqnarray*}  
e^{-\frac{n-1}{2}t} \| v_{\psi _0}  (x, \phi (t)) \| _{H_{(s)}}
& \leq  &
C  e^{-\frac{n-1}{2}t}\|\psi _0 \|_{H_{(s)}}  \quad \mbox{\rm for all} \,\,  t >0\,.
\end{eqnarray*}
where $\phi (t):= 1-e^{-t}$. 
For the second  term of (\ref{2.9}) we obtain
\begin{eqnarray*} 
&  &
e^{-\frac{n}{2}t}\int_{ 0}^{1} \|  v_{\psi _0}  (x, \phi (t)s) \| _{H_{(s)}} \big|2  K_0(\phi (t)s,t;M)+ nK_1(\phi (t)s,t;M)\big|\phi (t)\,  ds\\
& \leq &
 \|\psi _0 \|_{H_{(s)}}e^{-\frac{n}{2}t} \int_{ 0}^{1} \left( \big|2  K_0(\phi (t)s,t;M)\big|+ n\big|K_1(\phi (t)s,t;M)\big|\right) \phi (t)\,  ds \,.
\end{eqnarray*}
We have to estimate the following two integrals of the last inequality:
\begin{eqnarray*} 
&  & 
 \int_{ 0}^{1} 
\big|  K_i(\phi (t)s,t;M)\big|   \phi (t)\,  ds, \quad i=0,1 \,,
\end{eqnarray*}
where   $t>0$. To complete the estimate of the second term  of (\ref{2.9}) we are going to apply the next two lemmas
with $a=0$.
\begin{lemma}
\label{L2.3}
Let $  a>-1 $, $\Re M>0$,  and $\phi (t)= 1-e^{-t}$. Then 
\begin{eqnarray*} 
 \int_{ 0}^{1} \phi (t)^{a} s^{a}
\big|  K_1(\phi (t)s,t;M)\big|   \phi (t)\,  ds  
& \leq &
 C_M  e^{-at}(e^{t }-1)^{a+1} (e^{t }+1)^{ \Re M-1}  \quad \mbox{  for all}  \quad t>0\,.
\end{eqnarray*}  
In particular,
\begin{eqnarray*} 
 \int_{ 0}^{1} \phi (t)^{a} s^{a}
\big|  K_1(\phi (t)s,t;M)\big|   \phi (t)\,  ds  
& \leq &
 C_{M,a}   e^{\Re Mt}  \quad \mbox{  for large}  \quad t\,.
\end{eqnarray*}
\end{lemma}
\medskip

\noindent
{\bf Proof.} By the definition of the kernel $K_1$, we obtain
\begin{eqnarray*} 
&  &
 \int_{ 0}^{1} \phi (t)^{a} s^{a}
\big|  K_1(\phi (t)s,t;M)\big|   \phi (t)\,  ds   
   =    
 \int_{ 0}^{1-e^{-t}} r^{a} 
\big|  K_1(r,t;M)\big|    \,  dr \\
& \leq &  
   4 ^{-\Re M} e^{ \Re Mt }  \int_{ 0}^{1-e^{-t}} r^{a} 
  \big((1+e^{-t })^2 -   r  ^2\big)^{-\frac{1}{2}+\Re M    }\left|F\left(\frac{1}{2}-M   ,\frac{1}{2}-M  ;1; 
\frac{ ( 1-e^{-t })^2 -r^2 }{( 1+e^{-t })^2 -r^2 } \right)  \right|   \,  dr \\
& \leq &  
   4 ^{-\Re M} e^{ \Re Mt }  \int_{ 0}^{e^{t }-1} e^{t-2 \Re Mt }e^{-at}y^{a} 
  \big((e^{t }+1)^2 -   y  ^2\big)^{-\frac{1}{2}+\Re M    } \\
&  &
\hspace{2cm} \times \left| F\left(\frac{1}{2}-M   ,\frac{1}{2}-M  ;1; 
\frac{ ( e^{t }-1)^2 -y^2 }{( e^{t }+1)^2 -y^2 } \right)\right|     \, e^{-t } dy ,
\end{eqnarray*}
where the substitution  $e^tr=y$ has been used. Thus,
\begin{eqnarray*}   
 \int_{ 0}^{1} \phi (t)^{a} s^{a}
\big|  K_1(\phi (t)s,t;M)\big|   \phi (t)\,  ds  
& \leq &
   4 ^{-\Re M} e^{ -\Re Mt -at } \int_{ 0}^{e^{t }-1} y^{a} 
  \big((e^{t }+1)^2 -   y  ^2\big)^{-\frac{1}{2}+\Re M    } \\
  &  &
\hspace{2cm} \times  
\left| F\left(\frac{1}{2}-M   ,\frac{1}{2}-M  ;1; 
\frac{ ( e^{t }-1)^2 -y^2 }{( e^{t }+1)^2 -y^2 } \right)\right|     \,  dy \,.
\end{eqnarray*}
On the other hand, for $\Re M>0$ we have (see Section~\ref{A})
\begin{eqnarray*} 
&  &  
\left|F\left(\frac{1}{2}-M   ,\frac{1}{2}-M  ;1; \zeta  \right) \right|     \leq C_M \quad \mbox{\rm for all} \quad \zeta  \in [0,1)\,,
\end{eqnarray*}
where 
\begin{eqnarray*} 
\zeta := \frac{ ( e^{t }-1)^2 -y^2 }{( e^{t }+1)^2 -y^2 }  \in [0,1) \quad \mbox{\rm for all} \quad y \in [0, e^{t }-1] \quad \mbox{\rm and all} \quad t>0.
\end{eqnarray*}
Hence,
\begin{eqnarray*} 
 \int_{ 0}^{1} \phi (t)^{a} s^{a}
\big|  K_1(\phi (t)s,t;M)\big|   \phi (t)\,  ds 
& \leq &
 C_M  e^{-\Re Mt -at } \int_{ 0}^{e^{t }-1}  y^{ a} 
  \big((e^{t }+1)^2 -   y  ^2\big)^{-\frac{1}{2}+\Re M    }     \,   dy \,.
\end{eqnarray*}
If we denote $z:= e^t$, then for $M>0$ we have
\begin{eqnarray*} 
\int_{ 0}^{z-1}  y^{a} 
  \big((z+1)^2 -   y  ^2\big)^{-\frac{1}{2}+M    }     \,   dy  
& = &
\frac{1}{1+a}(z-1)^{1+a} (z+1)^{2 M-1} F\left(\frac{1+a}{2},\frac{1}{2}-M;\frac{3+a}{2};\frac{(z-1)^2}{(z+1)^2}\right) ,
 \end{eqnarray*}  
where $a>-1$ and $z\geq 1$. Hence, for $\Re M>0$ we have
\begin{eqnarray*}  
\int_{ 0}^{1} \phi (t)^{a} s^{a}
\big|  K_1(\phi (t)s,t;M)\big|   \phi (t)\,  ds   
&  \leq &
C_M  e^{-\Re Mt -at }(e^{t }-1)^{a+1} (e^{t }+1)^{2 \Re M-1}  \quad \mbox{\rm for all} \quad t>0\,.
\end{eqnarray*}
Thus the lemma is proved. \hfill $\square$

\begin{lemma}
Let $ a>-1$,  $\Re M>0$, and  $\phi (t)= 1-e^{-t}$. Then 
\begin{eqnarray*}  
 \int_{ 0}^{1} \phi (t)^{a} s^{a}
\big|  K_0(\phi (t)s,t;M)\big|   \phi (t)\,  ds  
& \leq &
  C_{M,a}  
(e^t-1)^{a+1} \times \cases{   e^ { -at} (e^t+1)^{ -\frac{1}{2}} \quad \mbox{\rm if} \quad \Re M<1/2\,, \cr
 e^ { (\Re M-a)t}  (e^t+1)^{ -1} \quad \mbox{\rm if} \quad \Re M > 1/2,} 
\end{eqnarray*} 
for all   $t>0$. In particular,
\begin{eqnarray*} 
 \int_{ 0}^{1} \phi (t)^{a} s^{a}
\big|  K_0(\phi (t)s,t;M)\big|   \phi (t)\,  ds  
& \leq & C_{M,a}  
\times \cases{   e^ {  \frac{1}{2}t}  \quad \mbox{\rm if} \quad \Re M<1/2\,, \cr
 e^ { \Re  M t}    \quad \mbox{\rm if} \quad \Re M> 1/2\,,}  
\end{eqnarray*}
for large  $t$.
\end{lemma}
\medskip

\noindent
{\bf Proof.} By substituting $K_0$ into integral, we obtain
\begin{eqnarray*} 
&  & 
 \int_{ 0}^{1} \phi (t)^{a} s^{a}
\big|  K_0(\phi (t)s,t;M)\big|   \phi (t)\,  ds \\
& \leq  &
4 ^ {-\Re M}  e^{ t \Re M} \int_{ 0}^{1-e^{-t}} r^{a}
   \big((1+e^{-t })^2 - r^2\big)^{  \Re M    } \frac{1}{ [(1-e^{ -t} )^2 -  r^2]\sqrt{(1+e^{-t } )^2 - r^2} }\\
&   &
\times  \Bigg|\Bigg[  \big(  e^{-t} -1 +M(e^{ -2t} -      1 -  r^2) \big) 
F \Big(\frac{1}{2}-M   ,\frac{1}{2}-M  ;1; \frac{ ( 1-e^{-t })^2 -r^2 }{( 1+e^{-t })^2 -r^2 }\Big) \\
&  &
\hspace{1cm}  +   \big( 1-e^{-2 t}+  r^2 \big)\Big( \frac{1}{2}+M\Big)
F \Big(-\frac{1}{2}-M   ,\frac{1}{2}-M  ;1; \frac{ ( 1-e^{-t })^2 -r^2 }{( 1+e^{-t })^2 -r^2 }\Big) \Bigg]\Bigg|    \,  d r \,.
\end{eqnarray*}
Now we   make the change $r=e^{-t}y$ in the last integral and obtain 
\begin{eqnarray*} 
&  & 
 \int_{ 0}^{1-e^{-t}} r^{a}
   \big((1+e^{-t })^2 - r^2\big)^{  \Re M    } \frac{1}{ [(1-e^{ -t} )^2 -  r^2]\sqrt{(1+e^{-t } )^2 - r^2} }\\
&   &
\times  \Bigg|\Bigg[  \big(  e^{-t} -1 +M(e^{ -2t} -      1 -  r^2) \big) 
F \Big(\frac{1}{2}-M   ,\frac{1}{2}-M  ;1; \frac{ ( 1-e^{-t })^2 -r^2 }{( 1+e^{-t })^2 -r^2 }\Big) \\
&  &
\hspace{1cm}  +   \big( 1-e^{-2 t}+  r^2 \big)\Big( \frac{1}{2}+M\Big)
F \Big(-\frac{1}{2}-M   ,\frac{1}{2}-M  ;1; \frac{ ( 1-e^{-t })^2 -r^2 }{( 1+e^{-t })^2 -r^2 }\Big) \Bigg]\Bigg|    \,  d r \\
& = &
 e^{ -2\Re Mt}  e^ {-at}\int_{ 0}^{e^{t}-1 }  y^{a}
    \big((e^{t }+1)^2 - y^2\big)^{  \Re M    } \frac{1}{ \big((e^{t }-1)^2 - y^2\big) \sqrt{(e^{t }+1)^2 - y^2} }\\
&   &
\times   \Bigg|\Bigg[ \big(  e^{t} -e^{2t} +M(1-e^{ 2t}  -  y^2) \big) 
F \Big(\frac{1}{2}-M   ,\frac{1}{2}-M  ;1; \frac{ (e^{t }- 1)^2 - y ^2 }{(e^{t }+ 1)^2 - y ^2 }\Big) \\
&  &
\hspace{1cm}  +   \big(e^{2t }-1+   y ^2 \big)\Big( \frac{1}{2}+M\Big)
F \Big(-\frac{1}{2}-M   ,\frac{1}{2}-M  ;1; \frac{ (e^{t }- 1)^2 - y ^2 }{(e^{t }+ 1)^2 - y ^2 }\Big) \Bigg]\Bigg|    \,  d  y  \,.
\end{eqnarray*}
Then we denote $z=e^t$ and derive
\begin{eqnarray*} 
&  & 
 \int_{ 0}^{1} \phi (t)^{a} s^{a}
\big|  K_0(\phi (t)s,t;M)\big|   \phi (t)\,  ds \\
& \leq &
   z^ {-(\Re M+ a)}\int_{ 0}^{z-1 }  y^{a}
    \big((z+1)^2 - y^2\big)^{  \Re M    } \frac{1}{ \big((z-1)^2 - y^2\big) \sqrt{(z+1)^2 - y^2} }\\
&   &
\times   \Bigg|\Bigg[ \big(  z -z^{2 } +M(1-z^{ 2 }  -  y^2) \big) 
F \Big(\frac{1}{2}-M   ,\frac{1}{2}-M  ;1; \frac{ (z- 1)^2 - y ^2 }{(z+ 1)^2 - y ^2 }\Big) \\
&  &
\hspace{1cm}  +   \big(z^{2 }-1+   y ^2 \big)\Big( \frac{1}{2}+M\Big)
F \Big(-\frac{1}{2}-M   ,\frac{1}{2}-M  ;1; \frac{ (z- 1)^2 - y ^2 }{(z+ 1)^2 - y ^2 }\Big) \Bigg]\Bigg|    \,  d  y  \,.
\end{eqnarray*}
To complete the proof of lemma we need the estimate given by the following proposition.

\begin{proposition}
If\, $a >-1 $ and $\Re M>0$, then 
\begin{eqnarray*} 
&  & 
 \int_{ 0}^{z-1 }  y^{a}
    \big((z+1)^2 - y^2\big)^{  \Re M    } \frac{1}{ \big((z-1)^2 - y^2\big) \sqrt{(z+1)^2 - y^2} }\\
&   &
\times   \Bigg|  \big(  z -z^{2 } +M(1-z^{ 2 }  -  y^2) \big) 
F \Big(\frac{1}{2}-M   ,\frac{1}{2}-M  ;1; \frac{ (z- 1)^2 - y ^2 }{(z+ 1)^2 - y ^2 }\Big) \\
&  &
\hspace{1cm}  +   \big(z^{2 }-1+   y ^2 \big)\Big( \frac{1}{2}+M\Big)
F \Big(-\frac{1}{2}-M   ,\frac{1}{2}-M  ;1; \frac{ (z- 1)^2 - y ^2 }{(z+ 1)^2 - y ^2 }\Big)  \Bigg|    \,  d  y  \\& \leq &
C_{M,n,p,q,s}(z-1)^{1+a}\times \cases{  (z+1)^{\Re M-\frac{1}{2}} \quad \mbox{\rm if} \quad \Re M<1/2\,,\cr
 (z+1)^{2\Re M-1} \quad \mbox{\rm if}  \quad \Re M> 1/2}\,.
\end{eqnarray*}
\end{proposition}
\medskip

\noindent
{\bf Proof.}   We follow the arguments have been used in the proof of Lemma~7.4~\cite{Yag_Galst_CMP}. For $\Re M>0$   both hypergeometric functions are bounded.  We divide the domain of integration into  two zones, 
\begin{eqnarray*} 
Z_1(\varepsilon, z) 
& := &
\left\{ (z,r) \,\Big|\, \frac{ (z-1)^2 -r^2   }{ (z+1)^2 -r^2 } \leq \varepsilon,\,\, 0 \leq r \leq z-1 \right\} \,,\\  
Z_2(\varepsilon, z) 
& := &
\left\{ (z,r) \,\Big|\, \varepsilon \leq  \frac{ (z-1)^2 -r^2   }{ (z+1)^2 -r^2 },\,\, 0 \leq r \leq z-1  \right\},
\end{eqnarray*}
and then  split the integral into  two parts,
\begin{eqnarray*}
\int_{ 0}^{z-1} \star  \, dr 
& = &
\int_{ (z,r) \in Z_1(\varepsilon, z)   }  \star  \, dr 
+ \int_{ (z,r) \in Z_2(\varepsilon, z)  }  \star   \, dr \,.
\end{eqnarray*} 
 In the first zone $Z_1(\varepsilon, z) $ we have
\begin{eqnarray*}
F\Big(\frac{1}{2}-M,\frac{1}{2}-M;1; \frac{ (z-1)^2 -y^2   }{ (z+1)^2 -y^2 }   \Big)  
 &  =  &
 1 + \left( \frac{1}{2}-M \right)^2\frac{ (z-1)^2 -y^2   }{ (z+1)^2 -y^2 }   
+ O\left(\left( \frac{ (z-1)^2 -y^2   }{ (z+1)^2 -y^2 }\right)^2\right),  \\  
F\Big(-\frac{1}{2}-M,\frac{1}{2}-M;1; \frac{ (z-1)^2 -y^2   }{ (z+1)^2 -y^2 }   \Big)  
 &  =  & 
 1 - \left( \frac{1}{4}-M^2 \right)\frac{ (z-1)^2 -y^2   }{ (z+1)^2 -y^2 }   
+ O\left(\left( \frac{ (z-1)^2 -y^2   }{ (z+1)^2 -y^2 }\right)^2\right)   . 
\end{eqnarray*} 
 We use the last formulas to estimate the term containing hypergeometric functions: 
\begin{eqnarray*} 
&  &   
\Bigg|  \big(  z -z^{2 } +M(1-z^{ 2 }  -  y^2) \big) 
F \Big(\frac{1}{2}-M   ,\frac{1}{2}-M  ;1; \frac{ (z- 1)^2 - y ^2 }{(z+ 1)^2 - y ^2 }\Big) \\
&  &
\hspace{1cm}  +   \big(z^{2 }-1+   y ^2 \big)\Big( \frac{1}{2}+M\Big)
F \Big(-\frac{1}{2}-M   ,\frac{1}{2}-M  ;1; \frac{ (z- 1)^2 - y ^2 }{(z+ 1)^2 - y ^2 }\Big)  \Bigg|   \\
& \leq  &
 \frac{1}{2}\big( (z-1)^2 - y ^2 \big)  \\
&  & 
 + \frac{1}{8} |2 M-1|\left|  y^2+ 2z(z-1)+ z^2-1  +2 M \left(3 y^2+2z(z-1)+ z^2-1 \right)  \right| \frac{ (z-1)^2 -y^2   }{ (z+1)^2 -y^2 }    \\
&  &
 +   \frac{1}{2}\big( (z-1)^2 - y ^2 \big) O\left(\left( \frac{ (z-1)^2 -y^2   }{ (z+1)^2 -y^2 }\right)^2\right) \,.
\end{eqnarray*} 
Hence, we have to consider the following two integrals, which can be easily  estimated,
\begin{eqnarray*}
A_1
&  :=  &
\int_{ (z,y) \in Z_1(\varepsilon, z)  }  y^{a}  \big((z+1)^2 - y^2\big)^{\Re M -\frac{1}{2}   }\,  dy   
  \,, \\
A_2
& := &
z^2  \int_{  (z,y) \in Z_1(\varepsilon, z) }  y^{a}  \big((z+1)^2 - y^2\big)^{ \Re  M -\frac{3}{2}    }  dy 
   ,
\end{eqnarray*} 
 for all  $z \in [1,\infty)$. Indeed, for $A_1$ we obtain 
 \begin{eqnarray*}
A_1
& \leq   &
\int_{ 0 }^{z-1}  y^{a}  \big((z+1)^2 - y^2\big)^{\Re M -\frac{1}{2}   }\,  dy \\
& = &
\frac{1}{1+a}(z-1)^{1+a} (z+1)^{2 \Re M-1} F\Big(\frac{1+a}{2},\frac{1}{2}-\Re M;\frac{3+a}{2};\frac{(z-1)^2}{(z+1)^2}\Big) \\
& \leq &
C_{M,n,p,q,s} (z-1)^{1+a} (z+1)^{2 \Re M-1}\,.
\end{eqnarray*} 
 Similarly, if $\Re M>0$,  then
  \begin{eqnarray}
  \label{1.9}
A_2
& \leq   &
z^2  \int_{ 0 }^{z-1}   y^{a}  \big((z+1)^2 - y^2\big)^{  \Re M -\frac{3}{2}    }  dy \nonumber \\
& = &
z^2\frac{1}{1+a}(z-1)^{1+a} (z+1)^{2 \Re M-3} F\Big(\frac{1+a}{2},\frac{3}{2}-\Re M;\frac{3+a}{2};\frac{(z-1)^2}{(z+1)^2}\Big) \,.
\end{eqnarray} 
\smallskip

Here and henceforth, if $A$ and $B$ are two non-negative quantities, we use $A \lesssim  B$ to denote the statement that $A\leq CB $ for some absolute constant $C>0$. 

It suffices to consider the case of real valued $M$. Then (\ref{F32}) and (\ref{1.9}) in the case of $M<1/2$ imply 
 \begin{eqnarray*}
A_2
& \lesssim   &
z^2\frac{1}{1+a}(z-1)^{1+a} (z+1)^{2 M-3}  z^{\frac{1}{2}-M}  
  \lesssim    
  (z-1)^{1+a} (z+1)^{  M-\frac{1}{2} }\,.
\end{eqnarray*} 
In the case of $M\geq  1/2$ due to (\ref{F32a}) 
we derive
 \begin{eqnarray*}
A_2
& \lesssim  &
z^2(z-1)^{1+a} (z+1)^{2 M-3}   
  \lesssim   
 (z-1)^{1+a} (z+1)^{2 M-1}  \,.
\end{eqnarray*} 
 Finally, for the integral over  the first zone $Z_1(\varepsilon, z) $ we  obtain
\begin{eqnarray*}
\int_{ (z,r) \in Z_1(\varepsilon, z)   }  \star\,  dr
& \lesssim  &
(z-1)^{1+a}\times \cases{  (z+1)^{\Re M-\frac{1}{2}} \quad \mbox{\rm if} \quad \Re M<1/2\,,\cr
 (z+1)^{2\Re M-1} \quad \mbox{\rm if}  \quad \Re M>  1/2\,.}
\end{eqnarray*} 
In the second zone we have 
\[
0< \varepsilon \leq  \frac{ (z-1)^2 -r^2   }{ (z+1)^2 -r^2 } < 1 \quad \mbox{\rm and}  \quad 
\frac{ 1  }{ (z-1)^2 -r^2 }  \leq  \frac{ 1   }{ \varepsilon[(z+1)^2 -r^2] }\,.
\]
Then, the hypergeometric functions for $\Re M>0$ obey the estimates
\[
\left| F\Big(-\frac{1}{2}-M,\frac{1}{2}-M;1; \zeta      \Big) \right|  \leq C \,\,  \mbox{\rm and}  \,\,   
\left| F\Big(\frac{1}{2}-M,\frac{1}{2}-M;1; \zeta   \Big) \right|  \leq C_M   \,\, \, \mbox{\rm for all}\,\,   \zeta  \in [\varepsilon ,1) .
\]
This allows us to  estimate  the integral over the second zone as follows:
\begin{eqnarray*} 
&  & 
 \int_{  (z,y) \in Z_2(\varepsilon, z) }  y^{a}
    \big((z+1)^2 - y^2\big)^{  \Re M    } \frac{1}{ \big((z-1)^2 - y^2\big) \sqrt{(z+1)^2 - y^2} }\\
&   &
\times   \Bigg|  \big(  z -z^{2 } +M(1-z^{ 2 }  -  y^2) \big) 
F \Big(\frac{1}{2}-M   ,\frac{1}{2}-M  ;1; \frac{ (z- 1)^2 - y ^2 }{(z+ 1)^2 - y ^2 }\Big) \\
&  &
\hspace{1cm}  +   \big(z^{2 }-1+   y ^2 \big)\Big( \frac{1}{2}+M\Big)
F \Big(-\frac{1}{2}-M   ,\frac{1}{2}-M  ;1; \frac{ (z- 1)^2 - y ^2 }{(z+ 1)^2 - y ^2 }\Big)  \Bigg|    \,  d  y  \\
& \lesssim &
z^2   \int_{  (z,y) \in Z_2(\varepsilon, z) }  y^{a}
    \big((z+1)^2 - y^2\big)^{  \Re M  -\frac{3}{2}  } \,  d  y  \\ 
& \lesssim  &
 z^2   \int_{ 0}^{z-1}  y^{a}
    \big((z+1)^2 - y^2\big)^{  \Re M  -\frac{3}{2}  } \,  d  y  \,. 
\end{eqnarray*}
Then we apply (\ref{1.9}) and Lemma~\ref{L1.6}:   
\begin{eqnarray*}  
z^2   \int_{  (z,y) \in Z_2(\varepsilon, z) }  y^{a}
    \big((z+1)^2 - y^2\big)^{ \Re  M  -\frac{3}{2}  } \,  d  y   
& \lesssim &
 (z-1)^{1+a}\times \cases{  (z+1)^{\Re M-\frac{1}{2}} \quad \mbox{\rm if} \quad \Re M<1/2\,,\cr
 (z+1)^{2\Re M-1} \quad \mbox{\rm if}  \quad \Re M> 1/2\,,}  
\end{eqnarray*} 
for all  $z \in [1,\infty)$.  Finally, for the integral over  the second zone $Z_2(\varepsilon, z) $ we obtain
\begin{eqnarray*}
\int_{ (z,r) \in Z_2(\varepsilon, z)   }  \star\,  dr
& \lesssim &
(z-1)^{1+a}\times \cases{  (z+1)^{\Re M-\frac{1}{2}} \quad \mbox{\rm if} \quad \Re M<1/2\,,\cr
 (z+1)^{2\Re M-1} \quad \mbox{\rm if}  \quad \Re M>  1/2\,.}
\end{eqnarray*}
The rest of the proof is a repetition of the above used arguments.   
Thus, the  proposition is proved.    \hfill $\square$
\medskip 

\noindent
{\bf Completion of the proof of Theorem~\ref{T13.2}.} Thus, if $\psi _1=0 $, then from (\ref{2.9}) we derive
\begin{eqnarray*}  
&  &
\| \psi  (x,t) \| _{H_{(s)}}  \\
&   \leq &  
e^{-\frac{n-1}{2}t} \|  v_{\psi _0}  (x, \phi (t)) \| _{H_{(s)}}   \\
&  &
+ \,  e^{-\frac{n}{2}t}\int_{ 0}^{1} \| v_{\psi _0}  (x, \phi (t)s)\| _{H_{(s)}}  \big|2  K_0(\phi (t)s,t;M)+ nK_1(\phi (t)s,t;M)\big|\phi (t)\,  ds\\
&   \lesssim  &  
C  e^{-\frac{n-1}{2}t}(1-e^{-t})^{a}\|\psi _0 \|_{H_{(s)}}   \\
&  &
+ \|\psi _0 \|_{H_{(s)}}e^{-\frac{n}{2}t} \int_{ 0}^{1}  \left( \big|2  K_0(\phi (t)s,t;M)\big|+ n\big|K_1(\phi (t)s,t;M)\big|\right) \phi (t)\,  ds \\
&  \lesssim &  
  e^{-\frac{n-1}{2}t}\|\psi _0 \|_{H_{(s)}}   \\
&  &
+  \|\psi _0 \|_{H_{(s)}}e^{-\frac{n}{2}t} \Bigg(  (e^{t }-1) (e^{t }+1)^{ \Re M-1} +    
(e^t-1)\times \cases{   (e^t+1)^{ -\frac{1}{2}} \quad \mbox{\rm if} \quad \Re M<1/2\,, \cr
 e^ {  \Re M  t}  (e^t+1)^{ -1} \quad \mbox{\rm if} \quad \Re M >  1/2}  \Bigg) \\
&  \lesssim &  
  e^{-\frac{n-1}{2}t}\|\psi _0 \|_{H_{(s)}}  \\
&  &
+  \|\psi _0 \|_{H_{(s)}}e^{-\frac{n}{2}t}(e^{t }-1) 
 \left(   (e^{t }+1)^{ \Re M-1}    +    
 \cases{    (e^t+1)^{ -\frac{1}{2}} \quad \mbox{\rm if} \quad \Re M<1/2\,, \cr
 e^ { \Re Mt}  (e^t+1)^{ -1} \quad \mbox{\rm if} \quad \Re M>  1/2}  \right)\,.  
\end{eqnarray*}
 In particular, for large $t$ we obtain
\begin{eqnarray*}  
\|  \psi  (x,t) \| _{H_{(s)}} 
&  \lesssim &
    \|\psi _0 \|_{H_{(s)}} e^{-\frac{n-1}{2}t} 
+   \|\psi _0 \|_{H_{(s)}}e^{-\frac{n}{2}t} e^{t }    \left(   e^{(\Re M-1)t }     +    
 \cases{   e^{ -\frac{1}{2}t} \quad \mbox{\rm if} \quad \Re M<1/2\,, \cr
 e^ { \Re Mt}  e^{-t} \quad \mbox{\rm if} \quad \Re M>  1/2}  \right)\\
& \lesssim  &
    \|\psi _0 \|_{H_{(s)}}\left(  e^{-\frac{n-1}{2}t} 
+ e^{(-\frac{n}{2}+1)t}     \left[   e^{(\Re M-1)t }     +    
 \cases{   e^{ -\frac{1}{2}t} \quad \mbox{\rm if} \quad \Re M<1/2\,, \cr
 e^ {\Re  Mt}  e^{-t} \quad \mbox{\rm if} \quad \Re M>  1/2}  \right]\right) \,.  
\end{eqnarray*}
In the case of $\psi _0=0 $ we have 
\begin{eqnarray*} 
\| \psi  (x,t) \| _{H_{(s)}}
 & =  & 
2e^{-\frac{n}{2}t}\|\int_{0}^1   v_{\psi _1 } (x, \phi (t) s) 
  K_1(\phi (t)s,t;M) \phi (t)\, ds\| _{H_{(s)}} \\
 & \leq   & 
2\|\psi _1 \|_{H_{(s)}}e^{-\frac{n}{2}t}\int_{0}^1 
 | K_1(\phi (t)s,t;M)| \phi (t)\, ds\,.
\end{eqnarray*}
Due to Lemma~\ref{L2.3}  we obtain
\begin{eqnarray*} 
\| \psi (x,t) \| _{H_{(s)}}
 & \lesssim  & 
 \|\psi _1 \|_{B^{s,q}_{p}}e^{-\frac{n}{2}t}  (e^{t }-1) (e^{t }+1)^{ \Re M-1}\,.
\end{eqnarray*}
Theorem is proved. \hfill $\square$

\subsection{$\bf  H_{(s)} ({\mathbb R}^n) -H_{(s)} ({\mathbb R}^n) $ Estimate for the time derivatives of energy solutions}

\begin{theorem} 
Consider the Cauchy problem
\begin{eqnarray*}
&  &
  \psi _{tt}  - e^{-2 t} A(x,\partial_x)  \psi  +   n   \psi  _t   +  m^2  \psi   =  0 ,
   \quad  \psi (x,0)=  \psi _0 (x)  , \,\,  \psi _t(x,0)= \psi _1 (x) ,
 \end{eqnarray*} 
where $ A(x,\partial_x) = \sum_{|\alpha |\leq 2} a_{\alpha }(x)\partial _x^\alpha $ is a  second order negative elliptic partial differential operator, $a_{\alpha } \in {\mathcal B}^\infty $,  and $m^2  \in{\mathbb R} $. Then, there is 
a number $C>0 $
such that
\[
    \|  \psi _t (t) \|_{H_{(s)}}  + e^{-t}  \| \psi  (t)\|_{H_{(s+1)}} 
 \leq   
 C\left(    \|  \psi (t)\|_{H_{(s)}}  +   e^{-\frac{n}{2}t} \| \psi  _1\|_{H_{(s)}}  +     e^{-\frac{n}{2}t} \|  \psi _0 \|_{H_{(s+1)}}     \right) \,\, \mbox{ for all}\,\, t>0. 
\]
\end{theorem}
\medskip

\noindent
{\bf Proof.}
The change of unknown function 
\[
 \psi =e^{-\frac{n}{2}}u,\qquad u =e^{\frac{n}{2}} \psi 
\]
simplifies the equation. Therefore, we 
consider the Cauchy problem  
\[ 
u_{tt} - e^{-2t}A(x,D) u -M^2 u= 0 ,\quad \quad u(x,0)= u_0 (x)\, , \quad u_t(x,0)=u_1 (x)\,, 
\]
with  smooth initial functions $ u_0 (x)$ and $ u_1 (x) $. Here $M^2=n^2/4-m^2 $. The equation leads to the following identity
\begin{eqnarray*} 
&  &
\frac{1}{2} \frac{d}{dt} \left\{ (u_{t},u_{t}) - e^{-2t}  ( A(x,\partial_x) u,  u )  -  M^2     (u,u ) \right\} - e^{-2t}( A(x,\partial_x) u,  u )=  0 \,.
\end{eqnarray*}
Since the operator $A(x,\partial_x) $ is negative, it follows 
\begin{eqnarray*} 
&  &
\frac{1}{2} \frac{d}{dt} \left\{ (u_{t},u_{t}) - e^{-2t}  ( A(x,\partial_x) u,  u )  -  M^2     (u,u ) \right\} \leq   0 \,.
\end{eqnarray*}
The integration in time gives 
\begin{eqnarray*} 
&  &
  (u_{t},u_{t}) - e^{-2t}  ( A(x,\partial_x) u,  u )  -  M^2     (u,u )  \leq   (u_1,u_1) 
- e^{-2t}  ( A(x,\partial_x) u_0 ,  u_0  )  -  M^2     (u_0 ,u_0  ) \,,
\end{eqnarray*}
and, consequently,
\begin{eqnarray*} 
&  &
 \|u_{t}(t)\|_{L^2}  + e^{-t}  \|u (t)\|_{H_{(1)}} \leq  C  (\|u(t)\|_{L^2} + \|u_1\|_{L^2} +  \|u_0 \|_{H_{(1)}} )  \,.
\end{eqnarray*}
Similarly, using a standard technique (see, e.g., \cite{Taylor}), one can obtain for every $ s \in {\mathbb R}$ the following estimate
\begin{eqnarray*} 
&  &
 \|u_{t}(t)\|_{H_{(s)}}   + e^{-t}  \|u (t)\|_{H_{(s+1)}} \leq  C_s (\|u(t)\|_{H_{(s)}} + \|u_{1} \|_{H_{(s)}} +  \|u_0 \|_{H_{(s+1)}} )   \,.
\end{eqnarray*}
Then for the function $ \psi $ we have 
\begin{eqnarray*}
    \|ne^{\frac{n}{2}t} \psi  (t) + 2e^{\frac{n}{2}t}\psi _t (t) \|_{H_{(s)}}    
& \leq  &
 C_s\left(    e^{\frac{n}{2}t}\|  \psi (t)\|_{H_{(s)}}  +   \| \psi _1\|_{H_{(s)}}  +     \|  \psi _0 \|_{H_{(s+1)}}     \right)  \,,
\end{eqnarray*}
while
\begin{eqnarray*}
 e^{-t}  \| \psi  (t)\|_{H_{(s+1)}} 
& \leq &
 C_s e^{-\frac{n}{2}t}(\|u(t)\|_{H_{(s)}} + \|u_1\|_{H_{(s)}} +  \|u_0 \|_{H_{(s+1)}} ) \\
& \leq &
 C_s \| \psi (t)\|_{H_{(s)}} + e^{-\frac{n}{2}t}(\|\psi _1\|_{H_{(s)}} +  \|\psi _0 \|_{H_{(s+1)}} )  \,.
\end{eqnarray*}
Thus, the  theorem is proved. \hfill $\square$
\begin{theorem}
\label{T1.7}
For   $ s \in {\mathbb R}$ the solution    $ \psi = \psi (x,t)$ of the Cauchy problem (\ref{1.7}) 
for $M^2 \in {\mathbb R}$ and  $\Re M \in (0, 1/2)$  
satisfies the following estimate 
\begin{eqnarray*}
    \|  \psi _t (t) \|_{H_{(s)}}   
& \leq  &
 C    e^{-\frac{n-1}{2}t}\left(   \| \psi _1\|_{H_{(s)}}  +      \|  \psi _0 \|_{H_{(s+1)}}     \right)  \,.
\end{eqnarray*}
If  $M^2 \in {\mathbb R}$ and  $\Re M >  \frac{1}{2}$ or $M=1/2 $, then 
\begin{eqnarray*}
    \|  \psi _t (t) \|_{H_{(s)}}   
& \leq  &
 C\left(    \|  \psi (t)\|_{H_{(s)}}  +   e^{-\frac{n}{2}t} \| \psi _1\|_{H_{(s)}}  +     e^{-\frac{n}{2}t} \| \psi _0 \|_{H_{(s+1)}}     \right) 
\end{eqnarray*}
and
\begin{eqnarray*}
    \|  \psi _t (t) \|_{H_{(s)}}   
& \leq  &
 C      e^{(\Re M-\frac{n}{2})t}      \left\{ \|\psi _0 \|_{H_{(s+1)}} +  \|\psi _1 \|_{H_{(s)}}\right\}   \,.
\end{eqnarray*}
\end{theorem}
\medskip

\noindent
{\bf Proof.} 
According to Theorem~\ref{T13.2}  if $\Re M \in (0, 1/2)$, then
\[ 
\|  \psi  (t) \|_{H_{(s)}}   
  \leq   
C   e^{   -\frac{n-1}{2} t}\Big\{   \| \psi _0   \|_{H_{(s+1)}} 
+ (1- e^{-t}) \|\psi _1  \|_{H_{(s)}}  
 \Big\}  \,.
\] 
Hence
\begin{eqnarray*}
    \|  \psi _t (t) \|_{H_{(s)}}   
& \lesssim  &
  \|  \psi (t)   \|_{H_{(s)}} +  Ce^{-\frac{n}{2}t}\left(     \|  \psi _1\|_{H_{(s)}}  +     \|   \psi _0 \|_{H_{(s+1)}}     \right) \\   
& \lesssim  &
   e^{   -\frac{n-1}{2} t}\Big(   \|\psi _0   \|_{H_{(s+1)}} 
+ (1- e^{-t}) \|\psi _1  \|_{H_{(s)}}  \Big)  
+  e^{-\frac{n}{2}t}\left(     \| \psi _1\|_{H_{(s)}}  +     \|  \psi _0 \|_{H_{(s+1)}}     \right) \,.
\end{eqnarray*}
Similarly we  can consider the case of  $\Re M >  \frac{1}{2}$ .  
The  theorem is proved. \hfill $\square$

\section{$\bf  H_{(s)} ({\mathbb R}^n) -H_{(s)} ({\mathbb R}^n) $  Estimates for  Equations with  Source}
\label{S2b}

We consider equations with $m \in {\mathbb C} $ and $\frac{n^2}{4}\geq m^2   $ although result can be similarly obtained for the case 
of large mass, that is, for $\frac{n^2}{4}\leq m^2 $. Recall $M:=(n^2/4-m^2)^{1/2} $. 
In fact, for the   case of large mass and  for the case of $ M \in [1/2,n/2)$   (that is $m^2 \in (0,(n^2-1)/4]$) one can consult \cite{NARWA}.
This is why in the present paper we focus on the case of $M \in (0,1/2)\cup (n/2, \infty) $ and some complex valued $M$. Thus, we are interested also in the Higgs boson equation,   in the massive scalar fields   as well as in the tachyons having $m^2 < 0$.
\begin{theorem} 
\label{T11.1} 
Let $ \psi = \psi (x,t)$ be a solution of the Cauchy problem  
\[ 
 \psi _{tt} +   n    \psi _t - e^{-2 t} A(x,\partial_x)  \psi   + m^2 \psi  =  f\,, \quad  \psi   (x,0)= 0\, , \quad  \psi _t(x,0)=0\,.
\] 
Then  the solution $ \psi = \psi (x,t)$ for $0< \Re M <1/2 $    
satisfies the following  estimate:
\[  
\|    \psi (x ,t) \| _{H_{(s)} }
  \leq  
Ce^{-\frac{n-1}{2}t}   
 \int_{ 0}^{t} e^{\frac{n-1}{2}b} 
   \|f(x,b)\|_{H_{(s)} }\,db   \quad \mbox{for all}\quad   t>0.
\]  
If either $\Re M> 1/2$ or $M=1/2$, then 
\begin{eqnarray*}
\| \psi (x ,t)  \| _{H_{(s)} }
&  \leq   &
 C_M  e^{(\Re M -\frac{n}{2})t}    \int_{ 0}^{t}  e^{ -(\Re M- \frac{n}{2})b }  \|f(x,b)\|_{H_{(s)} }  \,db     \quad \mbox{for all}\quad   t>0\,.   
\end{eqnarray*}

Moreover, for the derivative $\partial_t  \psi  (x,t) $ if $0< \Re  M <1/2 $,  
then the following estimate holds
\begin{eqnarray*}  
\| \partial_t \psi (x ,t) \| _{H_{(s)} }  
 & \leq   &  
 Ce^{-\frac{n-1}{2}t}\int_{ 0}^{t}
  e^{\frac{n+1}{2}b}  \|f(x,b)\|_{H_{(s)} }    \, db \quad \mbox{for all}\quad   t>0.
\end{eqnarray*}
If $  \Re  M> 3/2$ or $M=3/2$, then
\begin{eqnarray*}  
\| \partial_t  \psi (x ,t) \| _{H_{(s)} }  
 & \leq   & 
  C   e^{(   \Re M -\frac{n}{2})t}    \int_{ 0}^{t}  e^{ -(  \Re  M- \frac{n}{2})b  }     \|f(x,b)\|_{H_{(s)} }  \,db \\  
 &    &
+  C e^{-\frac{n-1}{2}t}\int_{ 0}^{t}
  e^{\frac{n+1}{2}b}  \|f(x,b)\|_{H_{(s)} }    \, db \quad \mbox{for all}\quad   t>0.
\end{eqnarray*}
\end{theorem}
\medskip

\noindent
{\bf Proof.}  The case of $M=1/2 $ is an evident consequence of  the representation (\ref{1.5}) and in the remaining part of the proof it is not discussed. 
 From (\ref{Phy}) we have
\begin{eqnarray*} 
 \psi (x,t) 
&  =  &
2  e^{-\frac{n}{2}t} \int_{ 0}^{t} db
  \int_{ 0}^{ e^{-b}- e^{-t}} dr  \,  e^{\frac{n}{2}b} v(x,r ;b)  4 ^{-M}  e^{ M(b+t) } \Big((e^{-t }+e^{-b})^2 - r^2\Big)^{-\frac{1}{2}+M    } \\
  &  &
\times F\Big(\frac{1}{2}-M   ,\frac{1}{2}-M  ;1; 
\frac{ ( e^{-b}-e^{-t })^2 -r^2 }{( e^{-b}+e^{-t })^2 -r^2 } \Big) \,. 
\end{eqnarray*}
According to (\ref{Blplq}) we can write 
\begin{eqnarray*} 
&  & 
\| v(x,r ;b)   \| _{H_{(s)} }  
\le C  \|f(x,b)\|_{H_{(s)} }  \quad \mbox{\rm for all} \quad  r \in [0,1]\,.
\end{eqnarray*}
 Hence,
\begin{eqnarray*}
\| \psi (x ,t)  \| _{H_{(s)} } 
&  \leq   &
2  e^{-\frac{n}{2}t} \int_{ 0}^{t} db
  \int_{ 0}^{ e^{-b}- e^{-t}} dr  \,  e^{\frac{n}{2}b}\| v(x ,r ;b) \| _{H_{(s)} }   4 ^{-\Re M}  e^{ \Re M(b+t) } \\
&  &
\times \Big((e^{-t }+e^{-b})^2 - r^2\Big)^{-\frac{1}{2}+\Re M    } 
\left| F\Big(\frac{1}{2}-M   ,\frac{1}{2}-M  ;1; 
\frac{ ( e^{-b}-e^{-t })^2 -r^2 }{( e^{-b}+e^{-t })^2 -r^2 } \Big) \right|\\
& \lesssim  &
  e^{ \Re M t  } e^{-\frac{n}{2}t} \int_{ 0}^{t} e^{\frac{n}{2}b} e^{ \Re M b  } \|f(x,b)\|_{H_{(s)} }  \,db 
 \int_{ 0}^{ e^{-b}- e^{-t}} \,  \\
  &  &
\times 
\Big((e^{-t }+e^{-b})^2 - r^2\Big)^{-\frac{1}{2}+\Re M    }  \left| F\Big(\frac{1}{2}-M   ,\frac{1}{2}-M  ;1; 
\frac{ ( e^{-b}-e^{-t })^2 -r^2 }{( e^{-b}+e^{-t })^2 -r^2 } \Big) \right| dr \,.  
\end{eqnarray*}
Following the outline of the proof of Lemma~\ref{L2.3} we set   $  r=ye^{-t}$ and  obtain 
\begin{eqnarray*} 
\| \psi (x ,t)  \|_{H_{(s)} }
&  \leq   &
 C_M e^{ -\Re M t  } e^{-\frac{n}{2}t}  \int_{ 0}^{t} e^{\frac{n}{2}b} e^{ \Re M b  }\|f(x,b)\|_{H_{(s)} } \,db   \\
  &  &
 \times 
\int_{ 0}^{ e^{t-b}- 1} \,   \Big(( e^{t-b}+1)^2 - y^2\Big)^{-\frac{1}{2}+\Re M    } \left| F\Big(\frac{1}{2}-M   ,\frac{1}{2}-M  ;1; 
\frac{ ( e^{t-b}-1)^2 -y^2 }{( e^{t-b}+1)^2 -y^2 } \Big) \right| dy \,.  
\end{eqnarray*}
In order to estimate the second integral we apply Lemma~\ref{LA.5} with $z=e^{t-b} >1$ and $a=0$.
Hence, the   estimate (\ref{2.13}) implies for $0< \Re M<1/2 $ the following estimate
\begin{eqnarray*}
\| \psi (x ,t)  \| _{H_{(s)} }
&  \leq   &
 C_M e^{ -\Re M t  } e^{-\frac{n}{2}t}   \int_{ 0}^{t} e^{\frac{n}{2}b} e^{ \Re M b  }\|f(x,b)\|_{H_{(s)} } \,db   \\
  &  &
 \times 
\int_{ 0}^{ e^{t-b}- 1} \,    \Big(( e^{t-b}+1)^2 - y^2\Big)^{-\frac{1}{2}+\Re M    } \\
&  &
\times \left| F\Big(\frac{1}{2}-M   ,\frac{1}{2}-M  ;1; 
\frac{ ( e^{t-b}-1)^2 -y^2 }{( e^{t-b}+1)^2 -y^2 } \Big) \right| dy \\
&  \lesssim   &
  e^{ -\Re M t  } e^{-\frac{n}{2}t}  \int_{ 0}^{t} e^{\frac{n}{2}b} e^{ \Re M b  }\|f(x,b)\|_{H_{(s)} }   (e^{t-b}-1)  e^{(t-b)(\Re M-\frac{1}{2}) }  \,db \\ 
& \lesssim   &
   e^{-\frac{n-1}{2}t}    \int_{ 0}^{t} e^{\frac{n-1}{2}b }  \|f(x,b)\|_{H_{(s)} }  
 \,db \,,  
\end{eqnarray*}
while for $ \Re M > 1/2$    the estimate (\ref{2.14}) implies
\begin{eqnarray*}
\| \psi (x ,t)  \| _{H_{(s)} }
&  \leq   &
 C_M e^{ -\Re M t  } e^{-\frac{n}{2}t} \int_{ 0}^{t} e^{\frac{n}{2}b} e^{ \Re M b  }\|f(x,b)\|_{H_{(s)} } \,db   \nonumber  \\
  &  &
 \times 
\int_{ 0}^{ e^{t-b}- 1} \,   \Big(( e^{t-b}+1)^2 - y^2\Big)^{-\frac{1}{2}+\Re M    }  \nonumber \\
&  &
\times \left| F\Big(\frac{1}{2}-M   ,\frac{1}{2}-M  ;1; 
\frac{ ( e^{t-b}-1)^2 -y^2 }{( e^{t-b}+1)^2 -y^2 } \Big) \right| dy  \nonumber \\
& \lesssim   &
  e^{ -\Re M t  } e^{-\frac{n}{2}t}   \int_{ 0}^{t} e^{\frac{n}{2}b} e^{ \Re M b  }\|f(x,b)\|_{H_{(s)} }    (e^{t-b}-1) (e^{t-b}+1)^{2 \Re M-1}\,db  \nonumber  \\
& \lesssim  &
  e^{ -\Re M t  } e^{-\frac{n}{2}t}   \int_{ 0}^{t} e^{\frac{n}{2}b} e^{ \Re M b  }\|f(x,b)\|_{H_{(s)} }    (e^{t-b}-1) e^{(2\Re M-1)t}e^{-(2 \Re M-1)b} \,db  \nonumber \\
& \lesssim   &
  e^{(\Re M -\frac{n}{2})t}    \int_{ 0}^{t}  e^{ -(\Re M- \frac{n}{2})b  }    \|f(x,b)\|_{H_{(s)} } \,db \,.     
\end{eqnarray*} 

In order to estimate the time  derivative  of the function $ \psi $ in (\ref{Phy}) we write
\begin{eqnarray}
\label{1.10c}  
\partial_t \psi (x,t) 
 & =  &
-\frac{n}{2} \psi (x,t) +  2   e^{-\frac{n}{2}t}\int_{ 0}^{t} db\,
  e^{\frac{n}{2}b} v(x,e^{-b}- e^{-t} ;b) E(e^{-b}- e^{-t},t; 0,b;M)   \nonumber \\
&  &
+ 2   e^{-\frac{n}{2}t}\int_{ 0}^{t} db\, \int_{ 0}^{ e^{-b}- e^{-t}} dr  \,  e^{\frac{n}{2}b} v(x,r ;b) \partial_t   E(r,t; 0,b;M) .   
\end{eqnarray}
Further
\begin{eqnarray*} 
E(e^{-b}- e^{-t},t; 0,b;M) 
& = & 
\frac{1}{2}   e^{\frac{1}{2}b+\frac{1}{2}t} 
\end{eqnarray*} 
implies
\begin{equation}
\label{1.10b} 
   2   e^{-\frac{n}{2}t}\int_{ 0}^{t} 
  e^{\frac{n}{2}b} v(x,e^{-b}- e^{-t} ;b) E(e^{-b}- e^{-t},t; 0,b;M) \, db 
  =    
 e^{-\frac{n-1}{2}t}\int_{ 0}^{t}
  e^{\frac{n+1}{2}b} v(x,e^{-b}- e^{-t} ;b) \, db \,.    
\end{equation}
Due to  (\ref{Blplq})   we have
\begin{eqnarray*}
\|  v(x,r ;b) \| _{H_{(s)} } \leq   C \|f(x,b)\|_{H_{(s)} } \quad \mbox{\rm for all}\,\,\, r \in(0, e^{-b}- e^{-t})\subseteq (0,1]\,.
 \end{eqnarray*}
 Hence (\ref{1.10b}) implies
\begin{eqnarray} 
&  &
 \|   2   e^{-\frac{n}{2}t}\int_{ 0}^{t} db\,
  e^{\frac{n}{2}b} v(x,e^{-b}- e^{-t} ;b) E(e^{-b}- e^{-t},t; 0,b;M) \| _{H_{(s)} } \nonumber \\
& \leq  & 
 e^{-\frac{n-1}{2}t}\int_{ 0}^{t}
  e^{\frac{n+1}{2}b} \|v(x,e^{-b}- e^{-t} ;b) \| _{H_{(s)} }   \, db  \nonumber \\
  \label{1.20}
& \leq  & 
 Ce^{-\frac{n-1}{2}t}\int_{ 0}^{t}
  e^{\frac{n+1}{2}b}  \|f(x,b)\|_{H_{(s)} }    \, db \,.
\end{eqnarray}   
 For the last term of the derivative $\partial_t \psi  $ in (\ref{1.10c}) we have
\begin{eqnarray*} 
&  & 
 \| 2   e^{-\frac{n}{2}t}\int_{ 0}^{t} db\, \int_{ 0}^{ e^{-b}- e^{-t}} dr  \,  e^{\frac{n}{2}b} v(x,r ;b) \partial_t   E(r,t; 0,b;M)\| _{H_{(s)} }\\
& \lesssim  & 
    e^{-\frac{n}{2}t}\int_{ 0}^{t} db\, e^{\frac{n}{2}b}   \|f(x,b)\|_{H_{(s)} }\int_{ 0}^{ e^{-b}- e^{-t}} dr  \,| \partial_t   E(r,t; 0,b;M) |\,.
\end{eqnarray*} 
\begin{proposition}
\label{P2.2}
If $\Re M>0 $, then 
\begin{eqnarray*} 
\int_{ 0}^{ e^{-b}- e^{-t}} | \partial_t   E(r,t; 0,b;M) |\,dr 
& \lesssim & 
\cases{  e^{-\frac{1}{2} t } e^{-  b}+ e^{(\Re M  -\frac{1}{2} )t } e^{- 3b}\quad \mbox{if}  \quad  \Re M< 1/2 ,\cr 
e^{  \Re M (t  -  b)}   \quad \mbox{if}  \quad  \Re M > 3/2,  }
\end{eqnarray*} 
for all $t\geq 0 $ and $b\geq 0 $ such that $b<t $.
\end{proposition}
\medskip

\noindent
{\bf Proof.} We have
\begin{eqnarray} 
&  &
\partial_t E(r,t; 0,b;M) 
\label{2.19} \\
& = & 
\left(  \partial_t 4 ^{-M}  e^{ M(b+t) } \Big((e^{-t }+e^{-b})^2 - r^2\Big)^{-\frac{1}{2}+M    }  \right) F\Big(\frac{1}{2}-M   ,\frac{1}{2}-M  ;1; 
\frac{ ( e^{-b}-e^{-t })^2 -r^2 }{( e^{-b}+e^{-t })^2 -r^2 } \Big)\nonumber  \\
&  &
+  4 ^{-M}  e^{ M(b+t) } \Big((e^{-t }+e^{-b})^2 - r^2\Big)^{-\frac{1}{2}+M    }  \partial_t F\Big(\frac{1}{2}-M   ,\frac{1}{2}-M  ;1; 
\frac{ ( e^{-b}-e^{-t })^2 -r^2 }{( e^{-b}+e^{-t })^2 -r^2 } \Big) \nonumber  \,.
\end{eqnarray}
First we consider the second term of the  equation (\ref{2.19}). If $\Re M> 1/2$ and $M\not= 1/2$, then
\begin{eqnarray*} 
&  &
\Bigg|  4 ^{-M}  e^{ M(b+t) } \Big((e^{-t }+e^{-b})^2 - r^2\Big)^{-\frac{1}{2}+M    }  \partial_t F\Big(\frac{1}{2}-M   ,\frac{1}{2}-M  ;1; 
\frac{ ( e^{-b}-e^{-t })^2 -r^2 }{( e^{-b}+e^{-t })^2 -r^2 } \Big)  \Bigg| \nonumber\\
& = &
\Bigg|  4 ^{-M}  e^{ M(b+t) } \Big((e^{-t }+e^{-b})^2 - r^2\Big)^{-\frac{1}{2}+M    } \\
&  &
\times \Bigg[-\frac{1}{\left(r^2 \left(-e^{2 (b+t)}\right)+2 e^{b+t}+e^{2 b}+e^{2 t}\right)^2}(1-2 M)^2 e^{b+t} \left(r^2 e^{2 (b+t)}+e^{2 b}-e^{2 t}\right) \\
&  &
\times F\left(\frac{3}{2}-M,\frac{3}{2}-M;2;\frac{\left(e^{-b}-e^{-t}\right)^2-r^2}{\left(e^{-b}+e^{-t}\right)^2-r^2}\right) \Bigg]\Bigg| \nonumber\\ 
&\lesssim  &
 e^{ ( \Re M-3)(b+t) } \Big((e^{-t }+e^{-b})^2 - r^2\Big)^{-\frac{5}{2}  + \Re M    }  \left| r^2 e^{2 (b+t)}+e^{2 b}-e^{2 t}\right | \\ 
& \lesssim  &
 e^{ ( \Re M-3)(b+t) } \Big((e^{-t }+e^{-b})^2 - r^2\Big)^{-\frac{5}{2}  + \Re M    }  \left(e^{2 t} -e^{2 b}\right ) 
\end{eqnarray*}
and for $\Re M > 3/2 $  we can use (\ref{1.18}) of Lemma~\ref{L1.9b} with $a=0$ to estimate the integral of the last term:
\begin{eqnarray*} 
&  & 
\int_{ 0}^{ e^{-b}- e^{-t}}  \Bigg|  e^{ M(b+t) } \Big((e^{-t }+e^{-b})^2 - r^2\Big)^{-\frac{1}{2}+M    }  \partial_t F\Big(\frac{1}{2}-M   ,\frac{1}{2}-M  ;1; 
\frac{ ( e^{-b}-e^{-t })^2 -r^2 }{( e^{-b}+e^{-t })^2 -r^2 } \Big)  \Bigg|  \, dr \nonumber \\
& \lesssim  & 
\left(e^{2 t} -e^{2 b}\right )e^{ (\Re M-3)(b+t) } \int_{ 0}^{ e^{-b}- e^{-t}}    \Big((e^{-t }+e^{-b})^2 - r^2\Big)^{-\frac{5}{2}  +\Re M    }   \, dr \nonumber \\
& \lesssim  &
  e^{ (\Re M-3)(b+t) }\left(e^{2 t} -e^{2 b}\right ) e^{-(2\Re M-4) (b+t)}    \left(e^{t}-e^{b}\right)    \left(e^{b}+e^{t}\right)^{2 \Re M-5}   \nonumber \\
& \lesssim  & 
   e^{  \Re M (t  -  b)}   \quad \mbox{if}  \quad  \Re M > 3/2    \,. 
\end{eqnarray*}

For the case of  $\Re M<1/2 $ we have 
\begin{eqnarray*} 
&  &
\Bigg|    e^{ M(b+t) } \Big((e^{-t }+e^{-b})^2 - r^2\Big)^{-\frac{1}{2}+M    }  \partial_t F\Big(\frac{1}{2}-M   ,\frac{1}{2}-M  ;1; 
\frac{ ( e^{-b}-e^{-t })^2 -r^2 }{( e^{-b}+e^{-t })^2 -r^2 } \Big)  \Bigg| \nonumber\\
& = &
\Bigg|   e^{ M(b+t) } \Big((e^{-t }+e^{-b})^2 - r^2\Big)^{-\frac{1}{2}+M    } \\
&  &
\times \Bigg[\frac{1}{\left(r^2 \left(-e^{2 (b+t)}\right)+2 e^{b+t}+e^{2 b}+e^{2 t}\right)^2} (1-2 M)^2 e^{b+t} \left(r^2 e^{2 (b+t)}+e^{2 b}-e^{2 t}\right) \, 4^{3-2 M}\\
&  &
\times  \left(\frac{ 4 e^{-b}e^{-t}}{ \left(e^{-b}+e^{-t}\right)^2-r^2}\right)^{M-1} F\left(M+\frac{1}{2},M+\frac{1}{2};2;\frac{\left(e^{-b}-e^{-t}\right)^2-r^2}{\left(e^{-b}+e^{-t}\right)^2-r^2}\right)\Bigg]\Bigg| \nonumber\\ 
& \lesssim  &
\Bigg|   e^{ M(b+t) } \Big((e^{-t }+e^{-b})^2 - r^2\Big)^{-\frac{1}{2}+M    } \\
&  &
\times \Bigg[\frac{1}{ ( e^{4 (b+t)} )\left(\left(e^{-b}+e^{-t}\right)^2-r^2 \right)^2}  e^{b+t} \left(r^2 e^{2 (b+t)}+e^{2 b}-e^{2 t}\right)  \left(\frac{  e^{-(b+t)} }{ \left(e^{-b}+e^{-t}\right)^2-r^2}\right)^{M-1} \Bigg]\Bigg| \nonumber\\
&  \lesssim  &
\left|   e^{ -2(b+t) } \Big((e^{-t }+e^{-b})^2 - r^2\Big)^{-\frac{3}{2}    }   \left(r^2 e^{2 (b+t)}+e^{2 b}-e^{2 t}\right) \right|\\
& \lesssim  &
\Big((e^{-t }+e^{-b})^2 - r^2\Big)^{-\frac{3}{2}    }   e^{-2 b} \qquad  \mbox{ \rm for all}\quad r \leq  e^{-b} - e^{-t }\,. 
\end{eqnarray*}
Thus, for the case of $\Re M<1/2$ we obtain 
\begin{eqnarray*}
&  &
\Bigg|    e^{ M(b+t) } \Big((e^{-t }+e^{-b})^2 - r^2\Big)^{-\frac{1}{2}+M    }  \partial_t F\Big(\frac{1}{2}-M   ,\frac{1}{2}-M  ;1; 
\frac{ ( e^{-b}-e^{-t })^2 -r^2 }{( e^{-b}+e^{-t })^2 -r^2 } \Big)  \Bigg| \nonumber\\
& \lesssim  &
\Big((e^{-t }+e^{-b})^2 - r^2\Big)^{-\frac{3}{2}    }   e^{-2 b} \qquad  \mbox{ \rm for all}\quad r \leq  e^{-b} - e^{-t }\,. 
\end{eqnarray*} 
Next we apply Lemma~\ref{L1.9}  with $a=0$ and derive
\begin{eqnarray*} 
&  &
\int_{ 0}^{ e^{-b}- e^{-t}}  \Bigg|  e^{ M(b+t) } \Big((e^{-t }+e^{-b})^2 - r^2\Big)^{-\frac{1}{2}+M    }   \partial_t F\Big(\frac{1}{2}-M   ,\frac{1}{2}-M  ;1; 
\frac{ ( e^{-b}-e^{-t })^2 -r^2 }{( e^{-b}+e^{-t })^2 -r^2 } \Big)  \Bigg|  \, dr \nonumber \\
&  & 
\lesssim    e^{-\frac{1}{2} t } e^{-  b}\qquad \mbox{if}  \quad   \Re M<1/2 \,. 
\end{eqnarray*}  
Now we consider the first term of the   equation (\ref{2.19}):
\begin{eqnarray*} 
&  &
\left| \left(  \partial_t 4 ^{-M}  e^{ M(b+t) } \Big((e^{-t }+e^{-b})^2 - r^2\Big)^{-\frac{1}{2}+M    }  \right) F\Big(\frac{1}{2}-M   ,\frac{1}{2}-M  ;1; 
\frac{ ( e^{-b}-e^{-t })^2 -r^2 }{( e^{-b}+e^{-t })^2 -r^2 } \Big) \right| \nonumber \\
& \lesssim  &
\Bigg| \Bigg(  M e^{M (b+t)} \left(\left(e^{-b}+e^{-t}\right)^2-r^2\right)^{M-\frac{1}{2}}\\
&  &
-2 \left(M-\frac{1}{2}\right) \left(e^{-b}+e^{-t}\right) e^{b M+(M-1) t} \left(\left(e^{-b}+e^{-t}\right)^2-r^2\right)^{M-\frac{3}{2}} \Bigg)\Bigg| \\
& = &
\Bigg| \frac{e^{M (b+t)} \left(\left(e^{-b}+e^{-t}\right)^2-r^2\right)^M \left(M r^2 e^{2 (b+t)}+e^{2 b} (M-1)-e^{b+t}-M e^{2 t}\right)}{\sqrt{ (e^{-t}+e^{-  b})^2-r^2} \left(r^2 \left(-e^{2 (b+t)}\right)+2 e^{b+t}+e^{2 b}+e^{2 t}\right)}\Bigg| \\
& \lesssim  &
    e^{(\Re M ) (b+t)} \left(\left(e^{-b}+e^{-t}\right)^2-r^2\right)^{\Re M-\frac{3}{2}}  \left|  e^{-2 t}  +e^{-(b+t)}+M e^{-2b}-  Me^{-2 t} -M r^2\right|  \\
& \lesssim  &
    e^{(\Re M ) (b+t)} \left(\left(e^{-b}+e^{-t}\right)^2-r^2\right)^{\Re M-\frac{3}{2}}  \left(e^{-b}+e^{-t}\right)^2  \\
& \lesssim  &
    e^{(\Re M ) (b+t)}\cases{  e^{-2b-\Re M b}  \left(\left(e^{-b}+e^{-t}\right)^2-r^2\right)^{ -\frac{3}{2}} \quad \mbox{if} \quad  \Re M<1/2 \,,\cr
 \left(e^{-b}+e^{-t}\right)^2 \left(\left(e^{-b}+e^{-t}\right)^2-r^2\right)^{\Re M-\frac{3}{2}} \quad \mbox{if}  \quad \Re M>  1/2\,. } 
\end{eqnarray*}
Finally
\begin{eqnarray*} 
&  &
\left| \left(  \partial_t 4 ^{-M}  e^{ M(b+t) } \Big((e^{-t }+e^{-b})^2 - r^2\Big)^{-\frac{1}{2}+M    }  \right) F\Big(\frac{1}{2}-M   ,\frac{1}{2}-M  ;1; 
\frac{ ( e^{-b}-e^{-t })^2 -r^2 }{( e^{-b}+e^{-t })^2 -r^2 } \Big) \right| \nonumber \\
& \lesssim &
    e^{ \Re M  (b+t)}\cases{  e^{-2b-\Re Mb}  \left(\left(e^{-b}+e^{-t}\right)^2-r^2\right)^{ -\frac{3}{2}} \quad \mbox{if} \quad  \Re M<1/2\,, \cr
 \left(e^{-b}+e^{-t}\right)^2 \left(\left(e^{-b}+e^{-t}\right)^2-r^2\right)^{\Re M-\frac{3}{2}}\quad \mbox{if}  \quad \Re M> 1/2 \,, } 
\end{eqnarray*}
and due to Lemmas~\ref{L1.9}, \ref{LA2b}  with $a=0$ 
\begin{eqnarray*} 
&  &
\int_0^{e^{-b}-e^{-t}} dr\,  \Bigg| \left(  \partial_t 4 ^{-M}  e^{ M(b+t) } \Big((e^{-t }+e^{-b})^2 - r^2\Big)^{-\frac{1}{2}+M    }  \right) \nonumber  \\
&  &
\hspace{3cm} \times F\Big(\frac{1}{2}-M   ,\frac{1}{2}-M  ;1; 
\frac{ ( e^{-b}-e^{-t })^2 -r^2 }{( e^{-b}+e^{-t })^2 -r^2 } \Big) \Bigg| \nonumber \\
& \lesssim   &
 \int_0^{e^{-b}-e^{-t}} dr  \,  e^{\Re M  (b+t)}\cases{  e^{-2b-\Re Mb}  \left(\left(e^{-b}+e^{-t}\right)^2-r^2\right)^{ -\frac{3}{2}} \quad \mbox{if}  \quad \Re  M<1/2\,, \cr
 \left(e^{-b}+e^{-t}\right)^2 \left(\left(e^{-b}+e^{-t}\right)^2-r^2\right)^{\Re M-\frac{3}{2}}\quad \mbox{if}  \quad \Re M>  1/2 \,,} \nonumber \\
& \lesssim   &
e^{\Re M  (b+t)}  \cases{ e^{-2b-\Re Mb}\int_0^{e^{-b}-e^{-t}}    \left(\left(e^{-b}+e^{-t}\right)^2-r^2\right)^{ -\frac{3}{2}} \,dr  \quad \mbox{if}  \quad \Re  M<1/2\,, \cr
 \left(e^{-b}+e^{-t}\right)^2 \int_0^{e^{-b}-e^{-t}}   \left(\left(e^{-b}+e^{-t}\right)^2-r^2\right)^{\Re M-\frac{3}{2}} \,dr \quad \mbox{if}  \quad \Re M>  1/2,} \nonumber \\
& \lesssim  &
e^{\Re M  (b+t)}  \cases{ e^{-2b-\Re Mb}  e^{-\frac{1}{2} t } e^{-b} \quad \mbox{if}  \quad \Re  M<1/2\,, \cr
 \left(e^{-b}+e^{-t}\right)^2  \left(e^t-e^b\right) \left(e^b+e^t\right)^{2 (\Re M )-3} e^{-(a+2 (\Re M )-2) (b+t)}\,\, \mbox{if}  \,\, \Re M>  1/2,} \nonumber \\
& \lesssim  &
  \cases{ e^{\Re M  t}    e^{-\frac{1}{2} t } e^{- 3b} \quad \mbox{if}  \quad \Re  M<1/2\,, \cr
  \left(e^{-b}+e^{-t}\right)^2  e^t\left(e^b+e^t\right)^{2  \Re M  -3} e^{-(a+ \Re M  -2) (b+t)}\,\, \mbox{if}  \,\, \Re M>1/2 ,} \nonumber   \\
& \lesssim   &
  \cases{ e^{\Re M  t}    e^{-\frac{1}{2} t } e^{- 3b} \quad \mbox{if}  \quad \Re  M<1/2\,, \cr
    e^{\Re M (t- b)}\quad \mbox{if}  \quad \Re M>  1/2 \,.} 
\end{eqnarray*}
Thus,
\begin{eqnarray*} 
&  &
\int_0^{e^{-b}-e^{-t}} dr\,  \Bigg| \left(  \partial_t 4 ^{-M}  e^{ M(b+t) } \Big((e^{-t }+e^{-b})^2 - r^2\Big)^{-\frac{1}{2}+M    }  \right) F\Big(\frac{1}{2}-M   ,\frac{1}{2}-M  ;1; 
\frac{ ( e^{-b}-e^{-t })^2 -r^2 }{( e^{-b}+e^{-t })^2 -r^2 } \Big) \Bigg| \nonumber \\
&   &    
\lesssim  \cases{ e^{\Re M  t}    e^{-\frac{1}{2} t } e^{- 3b} \quad \mbox{if}  \quad \Re  M<1/2\,, \cr
    e^{\Re M (t- b)}\quad \mbox{if}  \quad \Re M>  1/2 \,.} 
\end{eqnarray*}
The proposition is proved. \hfill $\square$

Then  estimating the norms  in the case of $\Re  M<1/2$ we   obtain   
\begin{eqnarray*} 
&  & 
 \| 2   e^{-\frac{n}{2}t}\int_{ 0}^{t} db\, \int_{ 0}^{ e^{-b}- e^{-t}} dr  \,  e^{\frac{n}{2}b} v(x,r ;b) \partial_t   E(r,t; 0,b;M)\| _{H_{(s)} }\\
& \lesssim   & 
    e^{-\frac{n}{2}t}\int_{ 0}^{t} db\, e^{\frac{n}{2}b}   \|f(x,b)\|_{H_{(s)} }\int_{ 0}^{ e^{-b}- e^{-t}} dr  \,  | \partial_t   E(r,t; 0,b;M) | \,.
\end{eqnarray*} 
Next we apply Proposition~\ref{P2.2} and obtain
\begin{equation}
\label{2.16}
 \| 2   e^{-\frac{n}{2}t}\int_{ 0}^{t} db\, \int_{ 0}^{ e^{-b}- e^{-t}} dr  \,  e^{\frac{n}{2}b} v(x,r ;b) \partial_t   E(r,t; 0,b;M)\| _{H_{(s)} } \nonumber \\
 \lesssim  
    e^{-\frac{n}{2}t}\int_{ 0}^{t} e^{\frac{n-2}{2}b } \|f(x,b)\|_{H_{(s)} }\,db .
\end{equation}
By collecting estimates  (\ref{1.20})  and (\ref{2.16}) we obtain the final estimate for $\| \partial_t  \psi (x ,t) \| _{H_{(s)} } $ in  the case of $\Re M<1/2$.
For the case of  $\Re M> 3/2$, due to Proposition~\ref{P2.2} and according to (\ref{Blplq})      we have
\begin{eqnarray*}
&  & 
 \| 2   e^{-\frac{n}{2}t}\int_{ 0}^{t} db\, \int_{ 0}^{ e^{-b}- e^{-t}} dr  \,  e^{\frac{n}{2}b} v(x,r ;b) \partial_t   E(r,t; 0,b;M)\| _{H_{(s)} }  \\
& \lesssim & 
    e^{-\frac{n}{2}t}\int_{ 0}^{t} db\, \int_{ 0}^{ e^{-b}- e^{-t}} dr  \,  e^{\frac{n}{2}b} \|v(x,r ;b) \| _{H_{(s)} }|\partial_t   E(r,t; 0,b;M)|  \nonumber \\
& \lesssim   &  
   e^{(\Re M -\frac{n}{2})t}\int_{ 0}^{t} e^{(\frac{n}{2}-\Re M )b}     \|f(x,b)\|_{H_{(s)} }
     \,   db  \nonumber   \,.
\end{eqnarray*}
The last estimate together with  (\ref{1.20})  implies the last statement of the theorem. \\

\noindent
{\bf The case of $M=3/2 $.} We consider the equation (\ref{1.10c}), 
where
\begin{eqnarray*} 
E\left(r,t; 0,b;\frac{3}{2} \right)
& = &
\frac{1}{4} e^{-\frac{b}{2}-\frac{t}{2}} \left(e^{2 b}+e^{2 t} -r^2 e^{2 (b+t)}  \right)\,,\\
E\left(e^{-b}-e^{-t},t; 0,b;\frac{3}{2}\right)
& = &
\frac{1}{2} e^{\frac{b+t}{2}}\,,\\
\partial_t E\left(r,t; 0,b;\frac{3}{2}\right)
& = &
 \frac{1}{8} e^{-\frac{b}{2}-\frac{t}{2}} \left(3 e^{2 t}-3 r^2 e^{2 (b+t)}-e^{2 b}\right)\,.
\end{eqnarray*}
Consequently, for the first term of (\ref{1.10c}) we have
\begin{eqnarray*}
\| \psi (x ,t)  \| _{H_{(s)} }
&  \leq   &
 C   e^{(\frac{3}{2} -\frac{n}{2})t}    \int_{ 0}^{t}  e^{ -(\frac{3}{2}- \frac{n}{2})b }  \|f(x,b)\|_{H_{(s)} }  \,db   \,,   
\end{eqnarray*} 
while for the second term the following estimate
\begin{eqnarray*} 
 &   &
\| 2   e^{-\frac{n}{2}t}\int_{ 0}^{t}
  e^{\frac{n}{2}b} v(x,e^{-b}- e^{-t} ;b) E(e^{-b}- e^{-t},t; 0,b;\frac{3}{2}) \, db\| _{H_{(s)} } \nonumber  \\
  & \lesssim  &
 e^{-\frac{n-1}{2}t}   \int_{ 0}^{t} 
  e^{\frac{n+1}{2}b}  \| v(x,e^{-b}- e^{-t} ;b)\| _{H_{(s)} }    \, db \nonumber \\
  & \lesssim  &  
 e^{-\frac{n-1}{2}t}   \int_{ 0}^{t} 
  e^{\frac{n+1}{2}b}\|f(x,b)\|_{H_{(s)} }  \, db  
\end{eqnarray*} 
holds. For the last term of (\ref{1.10c}) we obtain
\begin{eqnarray*}
&  &
\|2   e^{-\frac{n}{2}t}\int_{ 0}^{t} db\, \int_{ 0}^{ e^{-b}- e^{-t}} dr  \,  e^{\frac{n}{2}b} v(x,r ;b) \partial_t   E(r,t; 0,b;M) \| _{H_{(s)} }\\
  & \lesssim  & 
e^{-\frac{n}{2}t}\int_{ 0}^{t} \, db\, e^{\frac{n}{2}b}e^{-\frac{b}{2}-\frac{t}{2}} \int_{ 0}^{ e^{-b}- e^{-t}} dr  \,  \| v(x,r ;b) \| _{H_{(s)} } 
  \left(3 r^2 e^{2 (b+t)}+e^{2 b}-3 e^{2 t}\right)  \\
  & \lesssim & 
e^{-\frac{n}{2}t}\int_{ 0}^{t} \, db\, e^{\frac{n}{2}b}e^{-\frac{b}{2}-\frac{t}{2}}\|f(x,b)\|_{H_{(s)} } \int_{ 0}^{ e^{-b}- e^{-t}}   
 \left(3 r^2 e^{2 (b+t)}+e^{2 b}-3 e^{2 t}\right)  \,dr \\
  & \lesssim  & 
e^{-\frac{n}{2}t}\int_{ 0}^{t} \, db\, e^{\frac{n}{2}b}e^{-\frac{b}{2}-\frac{t}{2}}\|f(x,b)\|_{H_{(s)} } e^{- (b+t)} \left(e^t-e^b\right)  \left|-3  e^{b+t}+3 e^{2 b}-3 e^{2 t}\right|\\
  & \lesssim  & 
e^{-\frac{n+1}{2}t}\int_{ 0}^{t}  e^{\frac{n-1}{2}b}\|f(x,b)\|_{H_{(s)} } \left(e^{-b}-e^{-t}\right) \left(e^{t}+e^{ b}\right)^2  \, db\\
  &\lesssim  & 
e^{-\frac{n-3}{2}t}  \int_{ 0}^{t}  e^{\frac{n-3}{2}b}    
 \|f(x,b)\|_{H_{(s)} }  \, db\,.   
\end{eqnarray*}
The final estimate for this case is
\begin{eqnarray*} 
\| \partial_t \psi  (x,t)  \| _{H_{(s)} }
 & \lesssim   &
  e^{-\frac{n-3}{2}t}    \int_{ 0}^{t}  e^{ \frac{n-3}{2}b }  \|f(x,b)\|_{H_{(s)} }  \,db 
+   e^{-\frac{n-1}{2}t}   \int_{ 0}^{t} 
  e^{\frac{n+1}{2}b} \|f(x,b)\|_{H_{(s)} }    \, db  .   
\end{eqnarray*}
The theorem is proved. \hfill $\square$

\section{Global Existence. Small Data Solutions}
\label{S3}

\setcounter{equation}{0}

We are going to apply the Banach's fixed-point theorem.  
In order to estimate nonlinear terms we use the   Lipschitz condition (${\mathcal L}$).
First we consider the integral equation  (\ref{5.1}),  
where the function $ \psi _0 (x,t) \in C([0,\infty);L^q ({\mathbb R}^n))$ is given. 
Every solution to   the equation (\ref{NWE}) solves also the last integral equation with some function $   \psi _0 (x,t)$. We note here that any classical  
solution to the  equation (\ref{NWE}) solves also the integral equation
(\ref{5.1}) with some function $\psi _0(t,x)$, which is a classical  
solution to the Cauchy problem for the linear equation  (\ref{1.7}).
\medskip

The operator $G$ and the structure of the nonlinear term determine the solvability of the integral equation (\ref{5.1}). 
For the operator $G $ generated by the linear part of the equation (\ref{0.1}) with $m^2<0 $ the global solvability 
of the integral 
equation (\ref{5.1}) was studied in \cite{yagdjian_DCDS}. For the case of  $m^2<0 $ and the nonlinearity $F( \psi ) = c  | \psi  |^{\alpha +1} $, $c\not=0$,  the results of \cite{yagdjian_DCDS}
imply the nonexistence of the global solution even for arbitrary small   function $ \psi _0 (x,0) $ under some conditions on $n$, $\alpha $, and $M \in {\mathbb C}$.  
 
\medskip

Consider the Cauchy problem in the Sobolev space $ H_{(s)} ({\mathbb R}^n)$ with $ s > n/2$, which is an algebra.
In the next theorem operator ${\mathcal K} $ (\ref{0.5}) is generated by linear part of the equation (\ref{NWE}).

\begin{theorem}
\label{TIE}  
Assume that  $F(x,u)$  is Lipschitz continuous in the  space $H_{(s)} ({\mathbb R}^n)$, $ s > n/2$, $F(x,0)=0$, and also that $\alpha >0 $. 

\noindent 
$(i)$ Suppose that  $0< \Re M<1/2$ and  $\gamma \in [0,\frac{n-1}{2}] $.
  Then for every given function $  \psi _0(x ,t) \in X({\varepsilon ,s,\gamma }) $  such that  
\begin{eqnarray*} 
&  &
\sup_{t \in [0,\infty)}  e^{\gamma  t}  \| \psi _0(\cdot ,t) \|_{H_{(s)} ({\mathbb R}^n)}  < \varepsilon\,,  
\end{eqnarray*}
and for sufficiently small  $\varepsilon $,   \, the integral equation (\ref{5.1}) has a unique solution \, $  \psi  (x ,t) \in X({2\varepsilon ,s,\gamma })  $. For the solution one has  
\begin{eqnarray}
\label{3.2}  
\sup_{t \in [0,\infty)}  e^{\gamma  t}  \| \psi  (\cdot ,t) \|_{H_{(s)} ({\mathbb R}^n)}  < 2\varepsilon \,.
\end{eqnarray}

\noindent
$(ii)$ Suppose that  $\Re M \in [1/2,n/2)$. Then for every given function $  \psi _0(x ,t) \in X({\varepsilon ,s,\gamma_0 }) $,  
$\gamma_0> 0 $, such that 
\begin{eqnarray*} 
&  &
\sup_{t \in [0,\infty)}  e^{\gamma_0  t}  \| \psi _0(\cdot ,t) \|_{H_{(s)} ({\mathbb R}^n)}  < \varepsilon\,,  
\end{eqnarray*} 
 for every $\gamma  $ such that $\gamma \leq \gamma _0 $, $\gamma <(n/2-\Re M)/(\alpha +1) $, and for sufficiently small  $\varepsilon $,  the integral equation (\ref{5.1})   has a unique solution \, $  \psi (x ,t) \in X({2\varepsilon ,s,\gamma })  $. For the solution one has (\ref{3.2}). \\    

\noindent
$(iii)$ Suppose that  $\Re M>n/2$. Then for the function $ \psi _0(x ,t) \in X({\varepsilon ,s,\gamma }) $,   
$  \gamma <\frac{1}{\alpha +1} (\frac{n}{2} - \Re M)  $,    
a unique solution \, $ \psi (x ,t)   $ of the integral equation (\ref{5.1}) has   
the lifespan $T_{ls}$ that can be estimated from below by \[
   T_{ls}    
  \geq  
- \frac{ 1}{|\gamma | } \ln   \left( \sup_{\tau  \in [0,\infty)} e^{\gamma \tau }\| \psi _0(\cdot ,\tau ) \|_{H_{(s)}({\mathbb R}^n)}   \right)   
-  C(M,n,\alpha ,\gamma ) \,.
\] 
 with some constant $C(M,n,\alpha ,\gamma ) $.
\end{theorem}
\medskip

\noindent
{\bf Proof.} 
(i) Consider the mapping
\begin{eqnarray*} 
 S[ \psi ] (x,t)
& := &
 \psi _0(x,t) + 
G[ F( \cdot , \psi )] (x,t)    \,. 
\end{eqnarray*} 
We are going to prove that $S$ maps $X({R,s,\gamma })$ into itself and that $S$ is a contraction, provided that   $\varepsilon  $ and $R$ are sufficiently small. Consider the case of  $\Re M=\Re (\frac{n^2}{4}-m^2)^{1/2}<1/2$.  
Theorem~\ref{T11.1}    implies 
\begin{eqnarray*}  
\|  S[ \psi ]  (\cdot  ,t) \|_{ H_{(s)}  ({\mathbb R}^n)  }  
&  \leq   &
 \| \psi _0(\cdot ,t) \|_{H_{(s)} ({\mathbb R}^n)  }  +     
    \|   G[ F(  \psi )](\cdot ,t)  \|_{H_{(s)} ({\mathbb R}^n)  }   \\ 
&  \leq   &
 \| \psi _0(\cdot ,t) \|_{H_{(s)} ({\mathbb R}^n)  }  + C_M   e^{- \frac{n-1}{2}  t}    
\int_{ 0}^{t} e^{ \frac{n-1}{2}  b}      \|  F(\cdot ,  \psi )(\cdot ,b)  \|_{H_{(s)} ({\mathbb R}^n)  }  \,db   \,.
\end{eqnarray*}
Taking into account the Condition (${\mathcal L}$)  
we arrive at 
\begin{eqnarray*}  
\|   S[ \psi ] (\cdot  ,t) \|_{ H_{(s)}  ({\mathbb R}^n)  }  
&  \leq   &
  \| \psi _0(\cdot ,t) \|_{ H_{(s)}({\mathbb R}^n)  }  + C_M   e^{- \frac{n-1}{2}  t}    
\int_{ 0}^{t} e^{ \frac{n-1}{2}  b}    \| \psi (\cdot ,b ) \| _{H_{(s)}({\mathbb R}^n)}   ^{\alpha  +1}  \,db  \,.
\end{eqnarray*}
Then, for $\gamma \in {\mathbb R} $ we have  
\begin{eqnarray*} 
&  &
e^{\gamma t}\|   S[ \psi ] (x ,t) \|_{ H_{(s)}({\mathbb R}^n)  } \\
&  \leq   &
 e^{\gamma   t}\| \psi _0(\cdot ,t) \|_{H_{(s)}({\mathbb R}^n)  }  +  C_M   e^{ \gamma  t- \frac{n-1}{2}  t}    
\int_{ 0}^{t} e^{ \frac{n-1}{2}  b}  e^{ - \gamma (\alpha +1)b} \left( e^{ \gamma  b }     \|  \psi (\cdot ,b ) \| _{H_{(s)}({\mathbb R}^n)}   \right)^{\alpha  +1}  \,db \\
&  \leq   &
 e^{\gamma  t}\| \psi _0(\cdot ,t) \|_{H_{(s)}({\mathbb R}^n)  } 
 +  C_M       
 \left( \sup_{\tau  \in [0,\infty)} e^{\gamma  \tau  }     \| \psi (\cdot ,\tau  ) \| _{H_{(s)}({\mathbb R}^n)}   \right)^{\alpha  +1}  e^{ \gamma  t- \frac{n-1}{2}  t} \int_{ 0}^{t} e^{ \frac{n-1}{2}  b}  e^{ - \gamma (\alpha +1)b} \,db  .
\end{eqnarray*}
For $   \gamma \in [0,   \frac{n-1}{2}]  $ and $\alpha >0 $, the following function is bounded  
\begin{eqnarray} 
\label{Int}
&  &
 e^{ \gamma  t- \frac{n-1}{2}  t} \int_{ 0}^{t} e^{ \frac{n-1}{2}  b}  e^{ - \gamma (\alpha +1)b} \,db \leq C  \quad \mbox{\rm for all}\quad  t \in [0,\infty)\,.
\end{eqnarray}
Consequently, 
\begin{eqnarray*} 
&  &
\sup_{t  \in [0,\infty)}e^{\gamma t}\|   S[ \psi ] (x ,t) \|_{ H_{(s)}({\mathbb R}^n)  } \\
&  \leq   &
  \sup_{t  \in [0,\infty)}e^{\gamma t}\| \psi _0(\cdot ,t) \|_{H_{(s)}({\mathbb R}^n)  }  +  C_M      
 \left( \sup_{t  \in [0,\infty)} e^{\gamma   t  }     \| \psi (\cdot ,\tau  ) \| _{H_{(s)}({\mathbb R}^n)}   \right)^{\alpha  +1} .  \nonumber
\end{eqnarray*}
 Thus, the last inequality proves that the operator $S$ maps $X({R,s,\gamma})$ into itself if $\varepsilon  $ and $R$ are sufficiently small, namely, if
$\varepsilon  +C R^{\alpha +1} < R $. 
\medskip

It remains to prove that $S$ is a contraction mapping.
As a matter of fact, we just  apply the estimate (\ref{calM})     and get the contraction property from 
\[
 e^{\gamma t} \|S[ \psi ](\cdot,t) -  S[\widetilde{\psi }   ](\cdot,t) \|_{H_{(s)}({\mathbb R}^n) }
 \leq  
CR(t) ^{\alpha } d( \psi ,\widetilde{\psi } )\,,  
\]
where $\displaystyle   R(t):= \max\{ \sup_{0\leq \tau \leq t } e^{\gamma \tau  }  \| \psi (\cdot ,\tau ) \| _{H_{(s)}({\mathbb R}^n) } , \sup_{0\leq \tau \leq t } e^{\gamma \tau  }  \| \widetilde{\psi } (\cdot ,\tau ) \| _{H_{(s)}({\mathbb R}^n)  }\} \leq R$. 
Indeed,   we have
\begin{eqnarray*}
&  &
\| S[ \psi ](\cdot , t) -  S[\widetilde{\psi } ](\cdot , t) \|_{H_{(s)}({\mathbb R}^n) }  
 = 
\| G[ \,
( F(\cdot , \psi ) - F( \cdot ,\widetilde{\psi }  ) )  ](\cdot, t)
\|_{H_{(s)}({\mathbb R}^n) }  \\
& \le  &
C_M     e^{- \frac{n-1}{2}  t}  
\int_{ 0}^{t}  e^{ \frac{n-1}{2}  b} \|  ( F(\cdot , \psi  ) - F(\cdot ,\widetilde{\psi } ) )  (\cdot ,b)) \| _{H_{(s)}({\mathbb R}^n) } \,db \\
& \le  &
C_M      e^{- \frac{n-1}{2}  t}  
\int_{ 0}^{t}  e^{ \frac{n-1}{2}  b}\| \psi (\cdot , b) -\widetilde{\psi } (\cdot , b) \|_{H_{(s)}({\mathbb R}^n) }
\Big( \|  \psi  (\cdot ,  b) \|_{H_{(s)}({\mathbb R}^n) } ^\alpha  
+ \|\widetilde{\psi } (\cdot ,  b)  \|_{H_{(s)}({\mathbb R}^n) }^\alpha 
\Big) \,db    \,.
\end{eqnarray*}
Thus,  taking into account (\ref{Int}), the last estimate, and the definition of the metric $ d( \psi  ,\widetilde{\psi } ) $, we obtain
\begin{eqnarray*}
&  &
e^{\gamma   t }  \| S[ \psi ](\cdot , t) -  S[\widetilde{\psi } ](\cdot , t) \|_{H_{(s)}({\mathbb R}^n) }   \\ 
& \le  &
C_M     e^{\gamma   t }     e^{- \frac{n-1}{2}  t}  
\int_{ 0}^{t}  e^{ \frac{n-1}{2}  b} \| \psi (\cdot , b) -\widetilde{\psi } (\cdot , b) \|_{H_{(s)}({\mathbb R}^n) } 
\Big( \|  \psi  (\cdot ,  b) \|_{H_{(s)}({\mathbb R}^n) } ^\alpha  
+ \|\widetilde{\psi } (\cdot ,  b)  \|_{H_{(s)}({\mathbb R}^n) }^\alpha 
\Big) \,db   \\ 
& \le  &
C_M      e^{\gamma   t }     e^{- \frac{n-1}{2}  t}  
\int_{ 0}^{t}  e^{ \frac{n-1}{2}  b-\gamma (\alpha+1)b  } \Big( \max_{0 \le \tau  \leq b } e^{\gamma \tau }\| \psi (\cdot , \tau ) -\widetilde{\psi } (\cdot , \tau ) \|_{H_{(s)}({\mathbb R}^n) }  \Big)\\
&  &
\times 
\Big(  \Big(  \max_{0 \le \tau  \leq b }e^{\gamma  \tau }\| \psi (\cdot ,  \tau ) \|_{H_{(s)}({\mathbb R}^n) }  \Big)^\alpha  
+  \Big(  \max_{0 \le \tau  \leq b }e^{\gamma  \tau }\| \widetilde{\psi }  (\cdot ,  \tau ) \|_{H_{(s)}({\mathbb R}^n) }  \Big)^\alpha   
\Big) \,db   \\
& \le  &
C_{M,\alpha}   d( \psi ,\widetilde{\psi } ) R(t)^\alpha      
     e^{\gamma   t - \frac{n-1}{2}  t}  
\int_{ 0}^{t}  e^{ \frac{n-1}{2}  b-\gamma (\alpha+1)b}  \,db   \\
& \le &
 C_{M,\alpha}    d( \psi ,\widetilde{\psi } ) R(t)^\alpha   \,.
\end{eqnarray*}
Consequently, 
\begin{eqnarray*}
e^{\gamma  t }  \| S[ \psi ](\cdot , t) -  S[\widetilde{\psi } ](\cdot , t) \|_{H_{(s)}({\mathbb R}^n) }  
& \le  &
 C_{M,\alpha} \delta^{-1} R(t)^\alpha   d( \psi ,\widetilde{\psi } )  \,.
\end{eqnarray*}
Then we choose $\varepsilon $ and $ R$ such that $C_{M,\alpha}   R ^\alpha <1 $.  Banach's fixed point theorem completes the proof of the case (i) of Theorem~\ref{TIE}.
\smallskip

\noindent
(ii) Consider now the case of $\Re M>1/2$. 
Theorem~\ref{T11.1}   implies 
\begin{eqnarray*}  
\|  S[ \psi ]  (\cdot  ,t) \|_{ H_{(s)}  ({\mathbb R}^n)  } 
&  \leq   &
 \| \psi _0(\cdot ,t) \|_{H_{(s)} ({\mathbb R}^n)  }  +     
    \|   G[ F( \psi  )](\cdot ,t)  \|_{H_{(s)} ({\mathbb R}^n)  }   \\ 
&  \leq   &
 \| \psi _0(\cdot ,t) \|_{H_{(s)} ({\mathbb R}^n)  }  +   C_M  e^{(\Re M -\frac{n}{2})t}    \int_{ 0}^{t}  e^{ -(\Re M- \frac{n}{2})b }  \| G[ F(  \psi   )](x,b)\|_{H_{(s)}}  \,db \,.
\end{eqnarray*}
Taking into account the Condition (${\mathcal L}$)  
we arrive at 
\begin{eqnarray*} 
&  &
\|   S[ \psi ] (\cdot  ,t) \|_{ H_{(s)}  ({\mathbb R}^n)  } \\
&  \leq   &
  \| \psi _0(\cdot ,t) \|_{ H_{(s)}({\mathbb R}^n)  }  + C_M   e^{(\Re M -\frac{n}{2})t}    \int_{ 0}^{t}  e^{ -(\Re M- \frac{n}{2})b }  \|  \psi (\cdot ,b ) \| _{H_{(s)}({\mathbb R}^n)}   ^{\alpha  +1}  \,db  \\
&  \leq   &
  \| \psi _0(\cdot ,t) \|_{ H_{(s)} ({\mathbb R}^n)  }  +  C_M  
 e^{(\Re M -\frac{n}{2})t}    \int_{ 0}^{t}  e^{ -(\Re M- \frac{n}{2})b }  e^{ - \gamma (\alpha +1)b} \left( e^{ \gamma  b }     \|  \psi (\cdot ,b ) \| _{H_{(s)}({\mathbb R}^n)}   \right)^{\alpha  +1}  \,db  \,  .  
\end{eqnarray*}
Then, for $  \psi _0(x,t)  \in X(R,s,\gamma _0)$ and $\gamma \geq 0 $ we have  
\begin{eqnarray*} 
&  &
e^{\gamma t}\|   S[ \psi ] (x ,t) \|_{ H_{(s)}({\mathbb R}^n)  } \\
&  \leq   &
 e^{\gamma   t}\| \psi _0(\cdot ,t) \|_{H_{(s)}({\mathbb R}^n)  } 
 +  C_M   e^{ \gamma  t+ (\Re M -\frac{n}{2})  t}    
\int_{ 0}^{t} e^{ -(\Re M -\frac{n}{2})  b}  e^{ - \gamma (\alpha +1)b} \left( e^{ \gamma  b }     \| \psi (\cdot ,b ) \| _{H_{(s)}({\mathbb R}^n)}   \right)^{\alpha  +1}  \,db \\
&  \leq   &
 e^{\gamma   t}\| \psi _0(\cdot ,t) \|_{H_{(s)}({\mathbb R}^n)  } \\
&  &
 +  C_M       
 \left( \sup_{\tau  \in [0,t)} e^{\gamma  \tau  }     \| \psi (\cdot ,\tau  ) \| _{H_{(s)}({\mathbb R}^n)}   \right)^{\alpha  +1}    e^{ \gamma  t+(\Re M -\frac{n}{2})  t} \int_{ 0}^{t} e^{ -(\Re M -\frac{n}{2})  b}  e^{ - \gamma (\alpha +1)b} \,db  .
\end{eqnarray*}
If $\gamma =\frac{1}{\alpha +1}(\frac{n}{2}- \Re M -\delta )>0$, $\gamma \leq \gamma _0 $, and $\delta >0 $, then
\begin{eqnarray*} 
&  &
e^{\gamma t}\|   S[ \psi ] (x ,t) \|_{ H_{(s)}({\mathbb R}^n)  } \\
&  \leq   &
 e^{\gamma   t}\| \psi _0(\cdot ,t) \|_{H_{(s)}({\mathbb R}^n)  }  +  C_M       
 \left( \sup_{\tau  \in [0,t)} e^{\gamma  \tau  }     \| \psi (\cdot ,\tau  ) \| _{H_{(s)}({\mathbb R}^n)}   \right)^{\alpha  +1}   e^{ \gamma  t+(\Re M -\frac{n}{2})  t} \int_{ 0}^{t} e^{ \delta  b}   \,db  \\
&  \leq   &
 e^{\gamma_0   t}\| \psi _0(\cdot ,t) \|_{H_{(s)}({\mathbb R}^n)  }  +  C_M e^{-\gamma \alpha t}\delta ^{-1}        
 \left( \sup_{\tau  \in [0,t)} e^{\gamma  \tau  }     \| \psi (\cdot ,\tau  ) \| _{H_{(s)}({\mathbb R}^n)}   \right)^{\alpha  +1}     .
\end{eqnarray*}
In  follows $ \psi  \in X(R,s,\gamma  ) $ provided that $ R$ and $\varepsilon  $ are sufficiently small. 
We skip the remaining part of the proof since it is similar to the case (i).
\smallskip

\noindent
(iii) Consider now the case of $\Re M\geq n/2>1/2$ and $ \psi _0 (x,t)\in X(R, s,\gamma ) $. 
Theorem~\ref{T11.1}  implies 
\begin{eqnarray*} 
\|  S[ \psi ]  (\cdot  ,t) \|_{ H_{(s)}  ({\mathbb R}^n)  } 
&  \leq   &
 \| \psi _0(\cdot ,t) \|_{H_{(s)} ({\mathbb R}^n)  }  +     
    \|   G[ F( \psi  )](\cdot ,t)  \|_{H_{(s)} ({\mathbb R}^n)  }   \\ 
&  \leq   &
 \| \psi _0(\cdot ,t) \|_{H_{(s)} ({\mathbb R}^n)  }  +   C_M  e^{(\Re M -\frac{n}{2})t}    \int_{ 0}^{t}  e^{ -(\Re M- \frac{n}{2})b }  \| G[ F(  \psi  )](x,b)\|_{H_{(s)}}  \,db \,.
\end{eqnarray*}
Taking into account the Condition (${\mathcal L}$)  
we arrive at 
\begin{eqnarray*} 
e^{\gamma  t}\|   S[ \psi ] (\cdot  ,t) \|_{ H_{(s)}  ({\mathbb R}^n)  } 
&  \leq   &
e^{\gamma  t}  \| \psi _0(\cdot ,t) \|_{ H_{(s)}({\mathbb R}^n)  }  + C_M   e^{(\gamma +\Re M -\frac{n}{2})t}    \int_{ 0}^{t}  e^{ -(\Re M- \frac{n}{2})b }  \|  \psi (\cdot ,b ) \| _{H_{(s)}({\mathbb R}^n)}   ^{\alpha  +1}  \,db  \\
&  \leq   &
e^{\gamma  t}  \| \psi _0(\cdot ,t) \|_{ H_{(s)} ({\mathbb R}^n)  } \\
&  &
 +  C_M  
  \left(  \max_{\tau \in[0,t]} e^{\gamma  \tau }  \| \psi (\cdot ,\tau  ) \| _{H_{(s)}({\mathbb R}^n)}   \right)^{\alpha  +1} e^{(\gamma +\Re M -\frac{n}{2})t}   \int_{ 0}^{t}  e^{ -(\gamma (\alpha +1)+\Re M- \frac{n}{2})b }    \,db  \\
&  \leq   &
e^{\gamma  t}  \| \psi _0(\cdot ,t) \|_{ H_{(s)} ({\mathbb R}^n)  } \\
&  &
 +  C_M  
  \left(  \max_{\tau \in[0,t]} e^{\gamma  \tau }  \|  \psi (\cdot ,\tau  ) \| _{H_{(s)}({\mathbb R}^n)}   \right)^{\alpha  +1}    \frac{   e^{ - \gamma  \alpha  t } -e^{(\gamma +\Re M -\frac{n}{2})t}  }{-(\gamma (\alpha +1)+\Re M- \frac{n}{2})}   \,  .  
\end{eqnarray*} 
If 
$  \gamma < \frac{1}{\alpha +1}(\frac{n}{2}  - \Re M) \leq 0 $, 
then 
 for the given $  \psi _0 (x,t) \in  X({T, s,\gamma})$ the lifespan of the solution $ \psi $ can be estimated from below. Set 
 \begin{equation}
 \label{3.5b}
 T_\varepsilon :=\inf \{ T\,:\,\max_{\tau  \in [0,T]}e^{\gamma  \tau }  \|  \psi (x ,\tau ) \|_{ H_{(s)}({\mathbb R}^n)  } \geq 2\varepsilon \}\,, \quad \varepsilon :=   \max_{\tau  \in [0,\infty)} e^{\gamma  \tau } \| \psi _0(\cdot ,\tau ) \|_{H_{(s)}({\mathbb R}^n)  }  \,.
 \end{equation}
 Then
\[  
2\varepsilon    \leq   
\varepsilon  +  C_M  
 e^{-\gamma \alpha T_\varepsilon}      
 \varepsilon  ^{\alpha  +1} 
 \] 
 implies
\[  
  T_\varepsilon    
   \geq   
-\frac{ 1}{|\gamma|  } \ln  \frac{\varepsilon   }{ C_M } \,. 
\] 
Thus, the theorem is proved.  \hfill $\square$
\medskip

\begin{remark}
By the arguments have been used in the proof of the last theorem it is easy to derive the  existence of local (in time) solution even for large  initial data.
\end{remark}

 \noindent{\bf Proof of Theorem~\ref{T0.1}. } (i) The case of $\Re M\in (0,1/2)$. 
For the function   $ \psi _0 (x,t) $, that is, for the solution of the Cauchy problem (\ref{1.7})  and for $s>\frac{n}{2}$,   according to Theorem~\ref{T13.2}
we have the estimate 
\[ 
\|  \psi _0 (x,t) \|_{ { H}_{(s)} ({\mathbb R}^n)}
  \leq   
C_{M,n,s}  e^{   -\frac{n-1}{2} t}\Big\{   \| \psi _0   \|_{ { H}_{(s)} ({\mathbb R}^n)}
+  \|\psi _1  \|_{ { H}_{(s)} ({\mathbb R}^n)} 
 \Big\} \,.
\]
For every $T>0$ we have $ \psi _0 (x,t) \in C([0,T]; { H}_{(s)} ({\mathbb R}^n)) \bigcap C^1([0,T]; { H}_{(s-1)} ({\mathbb R}^n)) $. According to Theorem~\ref{TIE}, for every initial functions $\psi _0 $ and $\psi _1 $ 
the function $ \psi _0 (x,t) $ belongs to the space  $X(R,s,\frac{n-1}{2} ) $, where the operator $S $ is a contraction. 

\noindent
(ii)  In the case of $\Re M \in[1/2,n/2)$ for the function   $ \psi _0 (x,t)$, that is, for the solution of the Cauchy problem (\ref{1.7})  and for $s>\frac{n}{2}$,   according to Theorem~\ref{T13.2}
we have the estimate
\begin{eqnarray}
\label{3.5}
\| \psi _0  (x,t) \| _{ { H}_{(s)} ({\mathbb R}^n)} 
&   \leq &
C e^{(\Re M-\frac{n}{2})t} \left\{ \|\psi _0 \|_{ { H}_{(s)} ({\mathbb R}^n)} +\|\psi _1 \|_{ { H}_{(s)} ({\mathbb R}^n)} \right\}\,.
\end{eqnarray}
According to Theorem~\ref{TIE}, for every initial functions $\psi _0 (x)$ and $\psi _1 (x)$ 
the function $ \psi _0 (x,t)$ belongs to the space  $X(R,s,\frac{n}{2}-\Re M ) $. 

\noindent
(iii)  If  $\Re M > n/2 $, then according to Theorem~\ref{T13.2} for the solution of (\ref{1.7})
we have the estimate (\ref{3.5}) and, consequently,  
$ \psi _0 (x,t) \in X(R,s,\gamma )$ with $\gamma = n/2 -\Re M <0$ for some $R>0$. 
On the other hand,
\begin{eqnarray*}
e^{\gamma  t}\|  \psi (\cdot  ,t) \|_{ H_{(s)}  ({\mathbb R}^n)  } 
&  \leq   &
e^{\gamma  t}  \| \psi _0(\cdot ,t) \|_{ H_{(s)} ({\mathbb R}^n)  } \\
&  &
 +  C_M  
  \left(  \max_{\tau \in[0,t]} e^{\gamma  \tau }  \| \psi (\cdot ,\tau  ) \| _{H_{(s)}({\mathbb R}^n)}   \right)^{\alpha  +1}    \frac{ e^{ - \gamma  \alpha  t }-1 }{- \gamma  \alpha  } \,. 
\end{eqnarray*}
Next we define (\ref{3.5b}). 
Then  
\begin{eqnarray*} 
2 \varepsilon  
&  \leq   &
\varepsilon 
 +  C_M  
 \varepsilon ^{\alpha  +1}    
\frac{ e^{ - \gamma  \alpha  T_\varepsilon  }-1 }{- \gamma  \alpha  }  
\end{eqnarray*}
implies \,
$
T_\varepsilon \geq  -\frac{1}{ \Re M  - \frac{n}{2} }\ln \left(    \varepsilon  \right)- C(M,n,\alpha )
$. Theorem is proved.
\hfill $\square$
\medskip

 \noindent{\bf Proof of Theorem~\ref{T0.2}. } First consider the case of   $\Re  M\in(0,1/2)$. According to Theorem~\ref{T0.1}  there is a global solution $ \psi (x,t)  \in X(R,s,\frac{n-1}{2} ) $. Then, Theorem~\ref{T1.7} implies $\partial_t  \psi _0 (x,t)  \in X(R,s,\frac{n-1}{2} ) $. In order to check that $\partial_t F(\cdot ,\Phi) \in X(R,s,\frac{n-1}{2} ) $ we apply the Condition~($\mathcal L$) and the property of the operator $G$ proved in Theorem~\ref{T11.1} :
 \begin{eqnarray*}  
\| \partial_t G[F(x, \psi  )] (x ,t) \| _{H_{(s)}}  
 & \leq   & 
C e^{-\frac{n-1}{2}t}\int_{ 0}^{t} 
  e^{\frac{n+1}{2}b}  \|F(x, \psi  )(x,b)\|_{H_{(s)}}\,db\\  
 & \leq   & 
Ce^{-\frac{n-1}{2}t}   
 \int_{ 0}^{t} e^{\frac{n+1}{2}b}   
  \|  \psi (x,b)\|^{\alpha +1}_{H_{(s)}}\,db\,.
\end{eqnarray*}
 Hence, 
 \begin{eqnarray*}  
e^{\gamma _1 t}\| \partial_t G[F(x, \psi )] (x ,t) \| _{H_{(s)}}  
 & \leq   & 
C e^{(\gamma _1 -\frac{n-1}{2})t} \int_{ 0}^{t}   e^{\frac{n+1}{2}b}  e^{-(\alpha +1) \frac{n-1}{2}b}  
\left(  e^{\frac{n-1}{2}b}   
  \| \psi (x,b)\|_{H_{(s)}}\right)^{\alpha +1} \,db\\ 
 & \leq   & 
C \left( \max_{0\leq \tau \leq t} e^{\frac{n-1}{2}b}   
  \|  \psi (x,b)\|_{H_{(s)}}\right)^{\alpha +1} e^{(\gamma _1 -\frac{n-1}{2})t}\int_{ 0}^{t}  e^{1-\frac{\alpha }{2}(n-1)b}  \,db  \,.  
\end{eqnarray*}
If $\alpha >2/(n-1) $ we can set $\gamma _1= \frac{n-1}{2}$ and derive 
 \begin{eqnarray*}  
e^{\frac{n-1}{2}t}\| \partial_t G[F(x, \psi )] (x ,t) \| _{H_{(s)}}  
 & \leq   & 
C \left( \max_{0\leq \tau \leq t} e^{\frac{n-1}{2}b}   
  \| \psi (x,b)\|_{H_{(s)}}\right)^{\alpha +1} \,.  
\end{eqnarray*}
Assume now that  $  \Re  M> 3/2$ or $M=3/2$, then Theorem~\ref{T1.7} implies  $\partial_t  \psi _0 (x,t)  \in X(R,s, \frac{n}{2}- \Re M ) $.
Furthermore, 
\begin{eqnarray*}  
\| \partial_t G[F(x, \psi )] (x ,t) \| _{H_{(s)}}   
 & \leq   & 
  C   e^{(   \Re M -\frac{n}{2})t}    \int_{ 0}^{t}  e^{ -(  \Re  M- \frac{n}{2})b  }     \|  \psi (x,b)\|^{\alpha +1}_{H_{(s)}}  \,db \\  
 &    &
+  C e^{-\frac{n-1}{2}t}\int_{ 0}^{t}
  e^{\frac{n+1}{2}b}  \|  \psi (x,b)\|^{\alpha +1}_{H_{(s)}}   \, db \,,
\end{eqnarray*}
and
\begin{eqnarray*} 
e^{\gamma_1 t} \| \partial_t G[F(x,\psi )] (x ,t) \| _{H_{(s)}}  
 & \leq   & 
  C  \left( \max_{0\leq \tau \leq t} e^{\gamma b}   
  \|\psi (x,b)\|_{H_{(s)}}\right)^{\alpha +1} \\  
 &    &
\Big( e^{( \gamma_1+  \Re M -\frac{n}{2})t}    \int_{ 0}^{t}  e^{  ( \frac{n}{2} - \Re  M)b   - \gamma (\alpha +1)b  }   \,db 
+   e^{\gamma_1-\frac{n-1}{2}t}\int_{ 0}^{t}
  e^{\frac{n+1}{2}b  - \gamma (\alpha +1)b  }     \, db\Big) .
\end{eqnarray*}
Here $\frac{n}{2} - \Re  M    - \gamma (\alpha +1) >0$ and $\frac{n+1}{2}   - \gamma (\alpha +1)>0 $
since $\gamma <  \frac{1}{\alpha +1} (\frac{n}{2} - \Re  M) $.  The last factor is bounded if $\gamma _1   \leq \gamma (\alpha +1)-1 $. We set $\gamma _1 =\gamma -1<\min \{\frac{n}{2} - \Re  M,\gamma (\alpha +1)-1 \}$. 
 The theorem is proved. \hfill $\square$

\section*{Acknowledgment}
This paper was supported and completed within the  project University of Texas Rio Grande Valley College of Sciences 2016-17  Research Enhancement Seed Grant.

\begin{appendix}
\newtheorem{lemmaA}[theorem]{Lemma A}

\renewcommand{\theequation}{\thesection.\arabic{equation}}
\setcounter{equation}{0}

\section{Appendix}
\label{A}

There is a 
formula (see 15.3.6 of Ch.15\cite{A-S}  and \cite{B-E}) 
that ties together points $z=0$ and $z=1$:
\begin{eqnarray} 
\label{15.3.6}
 F  \left( a,b;c;z  \right) 
& = &
\frac{\Gamma (c)\Gamma (c-a-b)}{\Gamma (c-a)\Gamma (c-b)}F  \left( a,b;a+b-c+1;1-z  \right)  \\
 &  &
 + (1-z)^{c-a-b}\frac{\Gamma (c)\Gamma (a+b-c)}{\Gamma (a)\Gamma (b)}F  \left( c-a,c-b;c-a-b+1;1-z  \right) , \,\, |\arg (1-z)| <\pi . \nonumber  
\end{eqnarray}
Here $a,b,c \in{\mathbb C}$. It follows   
\begin{eqnarray}
\label{A2}
F  ( a,b;c;1  ) = \frac{ \Gamma (c)\Gamma (c-a-b) }{\Gamma (c-a)\Gamma (c-b)} \quad \mbox{\rm if } \quad \Re (c-a-b)>0.
\end{eqnarray}
 Each term of the   formula (\ref{15.3.6}) has a pole when $c=a+b\pm k$, ($k=0,1,2,\ldots$); this case is covered
by 15.3.10 of Ch.15\cite{A-S} 
\begin{eqnarray}
\label{A-S_15.3.10}
\hspace{-0.5cm} F  \left( a,b;a+b;z  \right) 
& = &
\frac{\Gamma (a+b)}{\Gamma (a)\Gamma (b)}\sum _{n=0}^\infty \frac{(a)_n(b)_n}{(n!)^2}
\left[ 2\psi (n+1) - \psi (a+n)- \psi (b+n) - \ln (1-z) \right] (1-z)^{n},
 \\
 &  &
\hspace*{5cm} \quad |\arg (1-z) |<\pi , \quad |1-z|<1 \,.\nonumber  
\end{eqnarray}

\begin{lemma}
\label{L1.6}
If $a>-1$ and  $M \in {\mathbb C} $ satisfies either $\Re \, M>1/2 $ or  $\Re M=1/2 \,\&\,\Im M\not=0$,  then
\begin{eqnarray}
\label{F32a}
\lim_{z \to \infty}  \,  F \left(\frac{a+1}{2},\frac{3}{2}-M;\frac{a+3}{2};\frac{(z-1)^2}{(z+1)^2}\right)
= \frac{\Gamma (\frac{a+3}{2})\Gamma (M - \frac{ 1}{2} )}{\Gamma (\frac{a }{2} +M)}  \,.
\end{eqnarray}
If $a>-1$  and $M=1/2 $, then
\begin{eqnarray*}
\lim_{z \to \infty}  \,\frac{1}{\ln z}  F \left(\frac{a+1}{2},\frac{3}{2}-M;\frac{a+3}{2};\frac{(z-1)^2}{(z+1)^2}\right)
= \frac{1 + a}{2} \,.
\end{eqnarray*}
If   $a>-1$ and $\Re \,M< 1/2$,  then
\begin{eqnarray}
\label{F32}
\lim_{z \to \infty} z^{M-\frac{1}{2}} \,  F \left(\frac{a+1}{2},\frac{3}{2}-M;\frac{a+3}{2};\frac{(z-1)^2}{(z+1)^2}\right)
=  2^{2 M-1}\frac{1+a }{1-2 M }\,.
\end{eqnarray}
\end{lemma}
\medskip

\noindent
{\bf Proof.} 
The statement (\ref{F32a}) follows from (\ref{A2}). 
Now consider the case of $\Re M<1/2$. According to \cite[(29) Sec.2.1.5]{B-E},
\begin{eqnarray*} 
  F \left(\frac{a+1}{2},\frac{3}{2}-M;\frac{a+3}{2};x \right)
& = &
(1-x)^{-( \frac{ 1}{2} -M )}  F \left(1,\frac{a }{2} +M;\frac{a+3}{2};x \right)
\end{eqnarray*}
while
\begin{eqnarray*} 
&  &
(1-x)^{-( \frac{ 1}{2} -M )} = 
\left(  \frac{4z}{(z+1)^2} \right)^{-( \frac{ 1}{2} -M )} \,,
\end{eqnarray*}
(\ref{A2}), and $1+\frac{a }{2} +\Re M<\frac{a+3}{2} $ yield
\begin{eqnarray*} 
& &
\lim_{z \to \infty} z^{M-\frac{1}{2}} \,  F \left(\frac{a+1}{2},\frac{3}{2}-M;\frac{a+3}{2};\frac{(z-1)^2}{(z+1)^2}\right)\\
& = &
\lim_{z \to \infty} z^{M-\frac{1}{2}} \, \left(  \frac{4z}{(z+1)^2} \right)^{-( \frac{ 1}{2} -M )} F \left(1,\frac{a }{2} +M;\frac{a+3}{2};\frac{(z-1)^2}{(z+1)^2}\right)\\
& = &
   4 ^{-( \frac{ 1}{2} -M )} \lim_{z \to \infty} F \left(1,\frac{a }{2} +M;\frac{a+3}{2};\frac{(z-1)^2}{(z+1)^2}\right)\\
& = &
   4 ^{-( \frac{ 1}{2} -M )} \lim_{\zeta  \to 1} F \left(1,\frac{a }{2} +M;\frac{a+3}{2};\zeta \right)\\
& = &
 4 ^{-( \frac{ 1}{2} -M )} \frac{ \Gamma \left(\frac{a+3}{2}\right) \Gamma \left(\frac{1}{2} -M\right) }{\Gamma \left(\frac{a+3}{2}-1\right) \Gamma \left(\frac{3 }{2} -M \right)}
= 4 ^{-( \frac{ 1}{2} -M )}  \frac{ \frac{a+1}{2} }{ \left(\frac{ 1}{2} -M \right)  }\,.
\end{eqnarray*}
If $M=1/2$, then we apply (\ref{A-S_15.3.10})
to  $
   F \left(\frac{a+1}{2},1;\frac{a+3}{2};\frac{(z-1)^2}{(z+1)^2}\right)
$ 
with 
\[
\zeta = \frac{(z-1)^2}{(z+1)^2}\,,\qquad 1-\zeta = \frac{4z}{(z+1)^2}\,,\qquad 
\lim_{z \to \infty}\frac{\ln(1-\zeta)}{\ln z} = -1\,.
\]
Hence,
\begin{eqnarray*} 
F  \left( \frac{a+1}{2},1;\frac{a+1}{2}+1;\zeta   \right)
& = &
\frac{\Gamma (\frac{a+3}{2} )}{\Gamma (\frac{a+1}{2})}\sum _{n=0}^\infty \frac{(\frac{a+1}{2})_n(1)_n}{(n!)^2}
\left[ \psi (n+1) - \psi (\frac{a+1}{2}+n)  - \ln (1-\zeta ) \right] (1-\zeta )^{n}\\
& = &
\frac{\Gamma (\frac{a+3}{2} )}{\Gamma (\frac{a+1}{2}) } 
\left[ \psi (1) - \psi (\frac{a+1}{2}) - \ln (1-\zeta ) \right]\\
&   &
+ \frac{\Gamma (\frac{a+3}{2} )}{\Gamma (\frac{a+1}{2})}\sum _{n=1}^\infty \frac{(\frac{a+1}{2})_n(1)_n}{(n!)^2}
\left[ \psi (n+1) - \psi (\frac{a+1}{2}+n)  - \ln (1-\zeta ) \right] (1-\zeta )^{n}
\end{eqnarray*}
and
\begin{eqnarray*} 
&  &
\lim_{z \to \infty}  \,\frac{1}{\ln z}  F \left(\frac{a+1}{2},1;\frac{a+3}{2};\frac{(z-1)^2}{(z+1)^2}\right)\\
&  = & 
\lim_{z \to \infty}  \frac{1}{\ln z} \Bigg\{\frac{\Gamma (\frac{a+3}{2} )}{\Gamma (\frac{a+1}{2})}  
\left[ \psi (1) - \psi (\frac{a+1}{2}) - \ln (1-\zeta ) \right]\\
&   &
+ \frac{\Gamma (\frac{a+3}{2} )}{\Gamma (\frac{a+1}{2})}\sum _{n=1}^\infty \frac{(\frac{a+1}{2})_n(1)_n}{(n!)^2}
\left[ \psi (n+1) - \psi (\frac{a+1}{2}+n)  - \ln (1-\zeta ) \right] (1-\zeta )^{n}\Bigg\}\\
&  = & 
 \frac{\Gamma (\frac{a+3}{2} )}{\Gamma (\frac{a+1}{2})} \lim_{z \to \infty} \frac{1}{\ln z}  
\left[ \psi (1) - \psi (\frac{a+1}{2}) - \ln (1-\zeta ) \right]\\
&   &
+ \lim_{z \to \infty} \frac{1}{\ln z} \frac{\Gamma (\frac{a+3}{2} )}{\Gamma (\frac{a+1}{2})}\sum _{n=1}^\infty \frac{(\frac{a+1}{2})_n(1)_n}{(n!)^2}
\left[  \psi (n+1) - \psi (\frac{a+1}{2}+n)  - \ln (1-\zeta ) \right] (1-\zeta )^{n} \\
&  = & 
- \frac{a+1}{2}   \lim_{z \to \infty} \frac{ \ln (1-\zeta )}{\ln z}  \\
& = &
\frac{1 + a}{2} \,.
\end{eqnarray*}
The lemma is proved. \hfill $\square$
\medskip

\begin{lemma}
\label{L1.9}
For $a >-1$ and $t>b>0$ we have
\begin{eqnarray*}
\int_0^{e^{-b}-e^{-t }} r^a \Big((e^{-t }+e^{-b})^2 - r^2\Big)^{-\frac{3}{2}    }   e^{-2 b} \,dr  
& \leq  &
C e^{-\frac{1}{2} t } e^{- (a+1)b}\,. 
\end{eqnarray*}
\end{lemma}
\medskip

\noindent
{\bf Proof.} We have 
\begin{eqnarray*} 
&  &
\int_0^{z-1}y^a ((z+1)^2-y^2)^{-3/2}dy \\
& = &
\frac{1}{4 (a+1) z (z+1)^3}(z-1)^{a+1} \Bigg\{   (a+2) (z+1)^2 F\left(-\frac{1}{2},\frac{a+1}{2};\frac{a+3}{2};\frac{(z-1)^2}{(z+1)^2}\right) \nonumber \\
&  &
-(z (4 a+z+2)+1) F\left(\frac{1}{2},\frac{a+1}{2};\frac{a+3}{2};\frac{(z-1)^2}{(z+1)^2}\right) \Bigg\} \nonumber \\
& \leq  &
C  z^{a-1}    \nonumber 
\end{eqnarray*}
for $z>1$. Then with $z=e^{t-b}$ we obtain
\begin{eqnarray*}
\int_0^{e^{-b}-e^{-t }} r^a \Big((e^{-t }+e^{-b})^2 - r^2\Big)^{-\frac{3}{2}    }   e^{-2 b} \,dr  
& = &
 e^{-2 b}  e^{-at -t+ 3/2t }  \int_0^{z-1} y^a \Big((z+1 )^2 - y^2\Big)^{-\frac{3}{2}    }  \,dy  \\
& \leq  &
 e^{-2 b}  e^{-at -t+ 3/2t }  z^{a-1} \\
& \leq  &
C e^{-\frac{1}{2} t }   e^{- (a+1)b}\,. 
\end{eqnarray*}
Lemma is proved. \hfill $\square$ 

\begin{lemma}
\label{LA2b}
Assume that $a>-1$, $t>b>0$, $M\in {\mathbb C}$,  and $\Re M\geq 1/2\,\&\, M\not=1/2$. Then
\begin{eqnarray*}  
&  &
\left| \int_0^{e^{-b}-e^{-t }} r^a \Big((e^{-t }+e^{-b})^2 - r^2\Big)^{-\frac{3}{2}  +  M  } \,dr \right| \\
& \leq   &
C(1+ (t-b)^{1-\sgn |\Re M-1/2|})\left(e^t-e^b\right)^{a+1} \left(e^b+e^t\right)^{2  \Re M -3} e^{-(a+2  \Re M -2) (b+t)} \,, 
\end{eqnarray*}
\end{lemma}
\medskip

\noindent
{\bf Proof.}  We have with $r= e^{-t }y $ 
\[
\int_0^{e^{-b}-e^{-t }} r^a \Big((e^{-t }+e^{-b})^2 - r^2\Big)^{-\frac{3}{2}  +M      } \,dr 
  =  
e^{-(a+1)t- (-3  +2M    )t }\int_0^{e^{t-b}-1} y^a \Big((1+e^{t-b})^2 - y)^2\Big)^{-\frac{3}{2}  +M     } \,dy  \,.
\]
 If $z:=e^{t-b} >1$, then we can  evaluate the last integral as follows:    
\[ 
\int_0^{z-1} y^a \left((1+z)^2 - y)^2\right)^{-\frac{3}{2}  +M   } \,dy  \\ 
  =   
\frac{1}{a+1}(z-1)^{a+1} (z+1)^{2 M -3} F \left(\frac{a+1}{2},\frac{3}{2}-M;\frac{a+3}{2};\frac{(z-1)^2}{(z+1)^2}\right)\,.
\] 
Hence, 
\begin{eqnarray*} 
&  &
\left| \int_0^{e^{-b}-e^{-t }} r^a \Big((e^{-t }+e^{-b})^2 - r^2\Big)^{-\frac{3}{2}  +M   } \,dr \right| \\ 
& = &
\frac{1}{a+1}\left(e^t-e^b\right)^{a+1} \left(e^b+e^t\right)^{2 \Re M-3} e^{-(a+2 \Re  M -2) (b+t)}\left|  F \left(\frac{a+1}{2},\frac{3}{2}-M ;\frac{a+3}{2};\frac{(e^{t-b}-1)^2}{(e^{t-b}+1)^2}\right)\right| \\
& \leq   &
C(1+ (t-b)^{1-\sgn |\Re M-1/2|})\left(e^t-e^b\right)^{a+1} \left(e^b+e^t\right)^{2  \Re M -3} e^{-(a+2  \Re M -2) (b+t)} \,. 
\end{eqnarray*}
Lemma is proved. \hfill $\square$

We skip the proof of the next lemma.
\begin{lemma}
\label{L1.9b}
Assume that $a>-1$, $t>b>0$,  $M\in {\mathbb C}$, $\Re\,M\geq 3/2$ and $M\not=3/2 $. Then
\begin{eqnarray}
\label{1.18}
&  &
\left| \int_0^{e^{-b}-e^{-t }} r^a \Big((e^{-t }+e^{-b})^2 - r^2\Big)^{-\frac{5}{2}  +M    } \,dr \right| \\
&  \leq  & 
C (1+ (t-b)^{1-\sgn |\Re M-3/2|})e^{-(a+2\Re M-4) (b+t)}    \left(e^{t}-e^{b}\right)^{1+a}    \left(e^{b}+e^{t}\right)^{2 \Re M-5} \nonumber    \,. 
\end{eqnarray}
\end{lemma}

\begin{lemma} 
\label{LA.5}If $\Re M>0$, $z>1$ and $a> -1$, then
\begin{eqnarray*} 
&  &
\int_{ 0}^{ z- 1} \, y^{a}   \Big(( z+1)^2 - y^2\Big)^{-\frac{1}{2}+\Re M    } \left| F\Big(\frac{1}{2}-M   ,\frac{1}{2}-M  ;1; 
\frac{ (z-1)^2 -y^2 }{(z+1)^2 -y^2 } \Big) \right| dy \\  
& \lesssim   & 
\cases{ (z-1)^{1+a}  z^{\Re M-\frac{1}{2} } \quad \mbox{\rm if} \quad  0<\Re M<1/2,  \cr 
 (z-1)^{1+a} (z+1)^{2 \Re M-1} \quad \mbox{\rm if} \quad    \Re M>1/2  .}
\end{eqnarray*}
\end{lemma}
\medskip

\noindent
{\bf Proof.} Since $\Re M>0$, then we have 
\begin{eqnarray*} 
&  &
\int_{ 0}^{ z- 1} \, y^{a}   \Big(( z+1)^2 - y^2\Big)^{-\frac{1}{2}+\Re M    } \left| F\Big(\frac{1}{2}-M   ,\frac{1}{2}-M  ;1; 
\frac{ (z-1)^2 -y^2 }{(z+1)^2 -y^2 } \Big) \right| dy \\
&  \lesssim    &
 \int_{ 0}^{ z- 1} \, y^{a}   \Big(( z+1)^2 - y^2\Big)^{-\frac{1}{2}+\Re M    }  dy\\
&  =   &
C_M\frac{1}{1+a}(z-1)^{1+a} (z+1)^{-1+2 \Re M} F\left(\frac{1+a}{2},\frac{1}{2}-\Re M;\frac{3+a}{2};\frac{(z-1)^2}{(z+1)^2}\right) \,  .
\end{eqnarray*}
Now we use Lemma~\ref{L1.6} (for $\Re M $), that is, 
\begin{eqnarray*} 
\lim_{z \to \infty} z^{M-\frac{1}{2}} \,  F \left(\frac{a+1}{2},\frac{3}{2}-M;\frac{a+3}{2};\frac{(z-1)^2}{(z+1)^2}\right)
=\frac{\pi  (a+1) 4^{M-1} \sec (\pi  M)}{\Gamma \left(\frac{3}{2}-M\right) \Gamma \left(M+\frac{1}{2}\right)}
\end{eqnarray*}
and obtain for $0<\Re M<1/2 $ 
\begin{eqnarray} 
&  &
\int_{ 0}^{ z- 1} \, y^{a}   \Big(( z+1)^2 - y^2\Big)^{-\frac{1}{2}+\Re M    } \left| F\Big(\frac{1}{2}-M   ,\frac{1}{2}-M  ;1; 
\frac{ (z-1)^2 -y^2 }{(z+1)^2 -y^2 } \Big) \right| dy \nonumber \\
&  \leq   & 
C_M\frac{1}{1+a}(z-1)^{1+a} (z+1)^{-1+2 \Re M} z^{\frac{1}{2}-\Re M} \nonumber \\
\label{2.13}
&  \leq   & 
C_M (z-1)^{1+a}  z^{\Re M-\frac{1}{2} } \,  ,
\end{eqnarray} 
while for $\Re M>1/2$   we obtain
\begin{eqnarray} 
&  &
\int_{ 0}^{ z- 1} \, y^{a}   \Big(( z+1)^2 - y^2\Big)^{-\frac{1}{2}+\Re M    } \left| F\Big(\frac{1}{2}-M   ,\frac{1}{2}-M  ;1; 
\frac{ (z-1)^2 -y^2 }{(z+1)^2 -y^2 } \Big) \right| dy \nonumber \\
&  \leq   & 
C \int_{ 0}^{ z- 1} \, y^{a}   \Big(( z+1)^2 - y^2\Big)^{-\frac{1}{2}+\Re M    } dy \nonumber \\
& = &
\frac{1}{1+a}(z-1)^{1+a} (z+1)^{2 \Re M-1} \left| F\Big(\frac{1+a}{2},\frac{1}{2}-\Re M;\frac{3+a}{2};\frac{(z-1)^2}{(z+1)^2}\Big)\right| \nonumber \\
\label{2.14}
& \lesssim  &
 (z-1)^{1+a} (z+1)^{2 \Re M-1}\,.
\end{eqnarray}
Lemma is proved. \hfill $\square$

\end{appendix}

\end{document}